\documentclass{siamltex}
\usepackage{amsfonts,latexsym}
\usepackage{amsmath}
\usepackage{graphicx}
\usepackage{amsmath}
\usepackage{color}
\usepackage{ulem}
\usepackage[top=3cm,bottom=4cm,right=4cm,left=4cm]{geometry}
\usepackage{booktabs} 
\usepackage{textcomp}
\usepackage{hyperref}
\usepackage{array}
\newtheorem{remark}{Remark}

\newtheorem{theorem1}{Theorem}
\newcommand{\R}{\mathbb R}
\newcommand{\C}{\mathbb C}

\newcommand{\cM}{\mathcal{M}}

\newcommand{\mR}{\mathcal{R}}

\newcommand{\bV}{\mathbb{V}}

\newcommand{\cV}{{\cal V}}
\newcommand{\bK}{\mathbb {K}}

\newcommand{\bH}{\mathbb {H}}
\newcommand{\cH}{{\cal H}}
\newcommand{\bT}{\mathbb {T}}

\newcommand{\bC}{\mathbb {C}}
\newcommand{\bB}{\mathbb {B}}
\newcommand{\bM}{\mathbb {M}}
\newcommand{\bA}{\mathbb {A}}

\newcommand{\cC}{{\cal C}}
\newcommand{\cB}{{\cal B}}

\newcommand{\cA}{{\cal A}}
\newcommand{\bp}{{\bf p}}
\newcommand{\bv}{{\bf v}}
\newcommand{\cJ}{{\cal J}}
\usepackage{epstopdf}
\usepackage{algorithm,algorithmic}

\title{An efficient extended block Arnoldi algorithm for feedback stabilization of incompressible Navier-Stokes flow problems}
\author{M.A. Hamadi \footnotemark[2] \thanks{Universit\'e Mohammed VI polytechnique, Lot 660, Hay Moulay Rachid, Ben Guerir, 43150 Maroc; amine.hamadi@um6p.ma; ahmed.ratnani@um6p.ma } \and K. Jbilou\footnotemark[1] \thanks{Universit\'e du Littoral C\^ote d'Opale, 50 Rue F. Buisson, BP 699--62228 Calais cedex, France;  jbilou@univ-littoral.fr} \and  A. Ratnani\footnotemark[1]}

\begin{document}
\maketitle
\begin{abstract}
Navier-Stokes equations are  well known in  modelling of an incompressible Newtonian fluid, such as air or water. This system of equations is very complex due to the non-linearity term that characterizes it. After the linearization and the discretization parts, we get a descriptor system of index-2 described by a set of differential algebraic equations (DAEs). The two main parts we develop  through this paper are focused firstly on constructing an efficient algorithm based on a projection technique onto an extended block Krylov subspace, that appropriately allows us to construct a reduced system  of the original DAE system. Secondly, we solve a Linear Quadratic Regulator (LQR) problem based on a Riccati feedback approach. This approach uses  numerical solutions  of  large-scale  algebraic Riccati equations. To this end, we use the  extended Krylov subspace method that allows us to project the initial large matrix problem onto a low order one that is solved by some direct methods. These numerical solutions are used to obtain a feedback matrix that will be used to stabilize the original system. We conclude by providing some numerical results to confirm the performances of our proposed method compared to other known methods.
\end{abstract}


\section{Introduction}
Navier-Stokes equations (NSEs) are very important  in the physics of fluid mechanics. The existence and smoothness of solutions is not yet guaranteed, although these equations are still of  interest to engineers and scientists in many technical fields. One of the main reason  that makes the solution of NSEs not unique is the chaotically appearing turbulences due  to  a naturally existing instabilities. In fact, these turbulences cannot be computed or predicted either, that is why we need to seek for  stabilization techniques. The stabilization of incompressible flow problems described by  Navier-Stokes equations is at the heart of a wide range of engineering applications, since they require a stable and controlled velocity field, which is considered to be the basis for ongoing reaction or production processes. A bench of work  based on the theoretical setting has been established by several authors for the stabilization of two and three-dimensional Navier-Stokes equations using a feedback control;  see M. Badra \cite{Badra1,Badra2}, V. Barbu et al.,\cite{Barbu1,Barbu2}, A.V. Fursikov \cite{Furs} and J. P. Raymond et al., \cite{Ray1,Ray2,Ray3}. Other works have been performed for the stabilization of two-dimensional Navier-Stokes equations based on a numerical setting by solving large-scale Linear Quadratic Regulator (LQR) problem using a Riccati-feedback approach; see \cite{benner15,Bansch}. In \cite{benner15} Bansch et al., proposed a generalized low-rank Cholesky factor Newton method to stabilize a flow around a cylinder. The LQR approach that interests us is based on a finite dimensional matrix derived from the discretization of  the linearized Navier-Stokes equations around a steady state. After the discretization stage we get a descriptor index-2 system of differential algebraic equations (DAEs) of a high dimension. Another way to deal with this stabilization problem is to  choose an appropriate method that allows us to construct a reduced system to the one described by a set of DAEs and then we stabilize the reduced system instead of the original one. This approach is convenient since it is based on the treatment of lower dimensional systems that makes the computation feasible. In \cite{Uddin} the authors use a balanced truncation method to construct an efficient reduced system and  they solve the obtained  LQR problem associated the reduced system based on a Riccati feedback approach.  

Two main parts will be covered  in this paper.  The first one focuses  on describing an efficient method to reduce a large-scale descriptor index-2 system of differential algebraic equations, depicted from a spatial discretization of the linearized Navier-Stokes equations around a steady state. This method is based on a projection technique onto an extended block Krylov subspace, and it allows us to construct a reduced system that has nearly the same response characteristics. A bench of work has been done to build an effective reduced model, such as  projection techniques onto  suitable Krylov-based subspaces as the rational or extended-rational  block Krylov subspaces, see \cite{houda1,frangos,Grim,hamadi,heyouni,jbilou,Druskin-on-optimal}.  Another class of methods described in \cite{Antoulas3,stykel,moore} contains balanced truncation methods. Numerous model reduction methods have been explored for  Navier-Stokes equations using balanced truncation and proper orthogonal decomposition \cite{bennerNstokes,Uddin}. A balanced truncation model reduction method for the Ossen equations has been investigated in \cite{bt}. For large problems, Krylov subspace methods are more efficient in term of cpu-time and memory requirements which is not the case for the methods based on balanced truncation since they require solving two large Lyapunov matrix equations at each iteration of the process and also the computation of singular value decompositions. All these methods that we mentioned here work properly for a class of descriptor dynamical systems represented by a set of ordinary differential equations (ODEs). Unfortunately, this is not our case since the dynamical system that we are dealing with is represented by a set of DAEs and therefore these methods are not directly applicable. Hence, one needs a process  that ensures a transformation of DAEs into ODEs in an appropriate manner. This will result in a dense projector called the Leray projection and  to overcome this problem, we give a simplification on how to avoid this dense projection matrix while  performing  our process to get a reduced system.  The second part of the paper is devoted to  solving a derived  Linear Quadratic Regulator (LQR) problem using a Riccati feedback approach. The major issue that we have to deal with is to solve a large-scale algebraic Riccati equation \cite{benricca,heyouni0,simon}, which is the key to design a controller represented by a feedback matrix. Our aim is to stabilize  the  unstable system by using the constructed feedback matrix. We propose an extended block Arnoldi algorithm  with appropriate computational requirements. The LQR problem used here is associated to the ODE system that relies on Leray  projections appearing  after the transformation to a set of an ODEs. We will explain how to avoid an explicit use of the Leray projection while solving  the obtained LQR problem. 

The remainder of this paper is structured as follows. In Section \ref{sec2}, we describe the incompressible Navier-Stokes equations with the linearization around a steady state, and its descriptor index-2 system of differential algebraic equations that arise after a mixed finite element method. The derivation of the obtained ODE system is also explained. Section \ref{sec3} deals with the extended block Krylov subspace method that allows us to construct an appropriate reduced model to the ODE system by avoiding the dense projection matrix that appears after the transformation to  ODEs. In Section \ref{sec4}, a Riccati feedback approach is explained  and we show how to  solve the LQR problem associated with the  ODE system. This approach is based on solving a large-scale algebraic  Riccati equation using an extended block Krylov subspace method. In the last section, we provide some numerical experiments to show the effectiveness of the proposed approaches.

\section{Navier-Stokes equations (NSEs) : Linearization and Discretization}
\label{sec2}

Navier-Stokes equations for a viscous, incompressible Newtonian fluid in a bounded domain $\Omega \subset \mathbb{R}^2$ with boundary $\partial \Omega$ are given by 
\begin{align}
	\label{NSequa}
	\left\{
	\begin{array}{lll}
		\dfrac{\partial z}{\partial t}- \dfrac{1}{\text{Re}} \Delta z+(z \cdot \nabla )z + \nabla p &=&  f,\\
		\nabla \cdot z   &=& 0,
	\end{array}
	\right.
\end{align}
where for   $t \in [0,\infty)$ and $x=[x_1 \quad x_2]^T \in \Omega \subset \R^2$, the vector $z(t,x)=[z_1(t,x),z_2(t,x)] \in \R^2$ refers to the velocity field, $p(t,x)\in \R$ is the pressure field, $f$ is known as the forcing term and  $\text{Re}  \in \R^{+}$ is the Reynolds number. The operators $\Delta$, $\nabla$ and  $\nabla\cdot$ are defined as the Laplacien, the Gradient and Divergence operators, respectively. The convective term in our model is a non-linear operator defined as 
\begin{align*}
	(z \cdot \nabla)z = \begin{bmatrix}
		z_1 \dfrac{\partial z_1}{\partial x_1} + z_2 \dfrac{\partial z_1}{\partial x_2} \\ 
		z_1 \dfrac{\partial z_2}{\partial x_1} + z_2 \dfrac{\partial z_2}{\partial x_2}
	\end{bmatrix}.
\end{align*}
The boundary $\Gamma = \partial \Omega$ can be partitioned as follows
$$\Gamma = \Gamma_{in} \cup \Gamma_{out} \cup \Gamma_{wall} \cup \Gamma_{feed}.$$
We therefore impose the following boundary conditions on the respective parts of the boundary 
$$
z = \left\{
\begin{array}{lll}
	\phi_{feed} & on & \Gamma_{feed}, \\
	\phi_{in} & on & \Gamma_{in}, \\
	0 & on & \Gamma_{wall}.
\end{array}
\right.
$$
The condition given below called, the \textit{do-nothing} condition 
$$-\dfrac{1}{\text{Re}} \nabla z \, n + p n = 0 \quad \text{on} \quad \Gamma_{out},$$ where $n$ denotes the unit outer normal vector to $\Gamma_{out}$.

\begin{figure}[H]
	\centering
	\includegraphics[width=12cm]{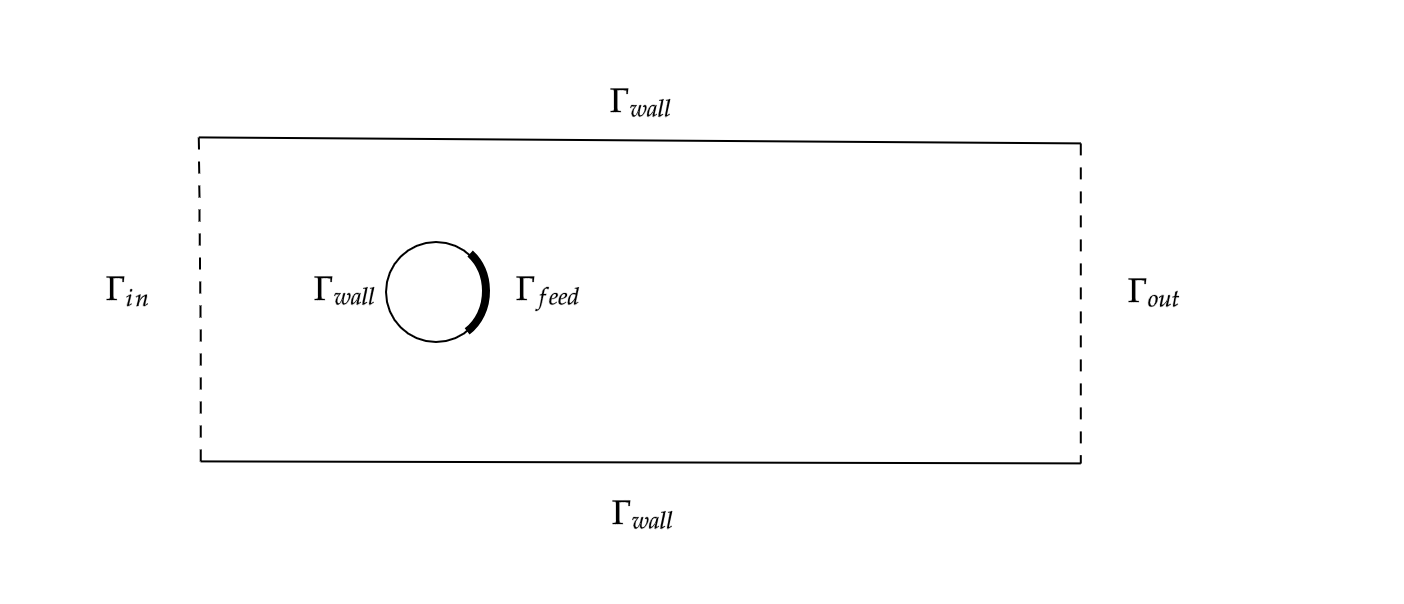}
	\caption{Domain $\Omega$ represented by a cylinder wake.}
\end{figure}
\noindent Navier-Stokes equations (NSEs) were derived independently by G.G. Stokes and C.L. Navier in the early 1800's. These equations describe the relationship between the velocity and the pressure of a moving fluid. NSEs represent the conservation of momentum. The fact that the convection term $(z \cdot \nabla) z$ is non-linear is what makes the NSEs complex. For incompressible flows, the second equation of the system (\ref{NSequa}) is called the continuity equation. In what follows, we present a linearization approach as it is described in \cite{benner15}.  
\subsection{Linearization}
We consider a stationary motion of an incompressible fluid described by the  velocity and pressure couple ($w_s(t,x),p_s(t,x)$)  that fulfills the stationary Navier-Stokes equations
\begin{eqnarray}
	-\dfrac{1}{\text{Re}} \Delta w_s +(w_s \cdot \nabla )w_s +\nabla p_s &=& f,\\ \nonumber
	\nabla \cdot w_s &=& 0.\nonumber
\end{eqnarray}
Here, the same boundary and initial conditions of the first equations are considered. The pair ($w_s,p_s$) depicts the desired stationary but possibly unstable solution of  system (\ref{NSequa}). \\
We define the following differences 
\begin{eqnarray*}
	v(t,x) &=& z(t,x)- w_s(x), \\
	\chi(t,x) &=& p(t,x)-p_s(t,x).
\end{eqnarray*}
Replacing in (\ref{NSequa}) and dropping the non-linear term, we obtain the following linearized Navies-Stokes equations
\begin{subequations}
	\begin{align}
		\dfrac{\partial v}{\partial t}-\dfrac{1}{\text{Re}} \Delta v +(w_s \cdot \nabla )v +(v \cdot \nabla )w_s  +\nabla \chi &= 0,\\
		\nabla \cdot v &= 0,
	\end{align}
	defined for $t\in [0,\infty)$ and $x \in \Omega \subset \mathbb{R}^2$ with Drichlet boundary conditions 
	\begin{eqnarray}
		v &=& 0 \quad \text{on} \quad \Gamma_{in} \cup \Gamma_{wall},  \\  
		v &=& \phi_{feed} \quad \text{on} \quad \Gamma_{feed},
	\end{eqnarray}
	a do-nothing condition is described as 
	$$-\dfrac{1}{\text{Re}} \nabla v \, n + \chi \, n = 0 \quad \text{on} \quad \Gamma_{out},$$
	and the initial condition 
	$$v(0,\cdot) = 0 \quad \text{in} \, \, \Omega.$$
\end{subequations}
$v$ is defined as perturbation of our flow field $z$ from the desired stationary flow field $w_s$. A \textit{zero output} for $t \to \infty$ implies that $v$  approximates $w_s$ for $t \to \infty$. As a consequence our flow field achieves the properties of the desired stationary flow field.
\subsection{The discrete equations}
The choice of an appropriate discretization technique depends on the specific governing equations used, for example (compressible or incompressible flow (our case), mesh type (structured, unstructured). The classical discretization techniques are finite difference,  finite element and finite volume. One of the known methods used to discretize instationary problems, is the method of lines  which is based on the replacement of the spatial derivatives in the PDE with algebraic approximations leading to a system of ODEs that approximates the original PDE. In this subsection, we briefly present the main properties of the discrete equations already established in \cite{benner15}. After using a mixed finite element method, we obtain a system of differential algebraic equations of the form 
\begin{subequations}
	\label{DAE1}
	\begin{align} 
		M \, \dfrac{d}{dt}\textbf{v}(t) &= A\textbf{v}(t)+G\textbf{p}(t)+\textbf{f}(t), \label{eq1} \\
		0 &= G^T\textbf{v}(t), \label{cond1}
	\end{align}
\end{subequations}
where 
\begin{enumerate}
	\item [] $\textbf{v}(t) \in \mathbb{R}^{n_v}$ : the nodal vector of the discretized velocity. \\
	\item[] $\textbf{p}(t) \in \mathbb{R}^{n_p}$ : the discretized pressure. \\
	\item[] $\textbf{f}(t) \in \mathbb{R}^{n_v}$ : the forcing term that contains the control. 
\end{enumerate}
\medskip
\noindent In what follows, we assume  that the forcing term $f(t)$ is given by
$$\textbf{f}(t) = B {\bf u}(t).$$
Moreover, the matrices $M=M^T \, \succ 0 \, \in \mathbb{R}^{n_v \times n_v}$ and $A \in \mathbb{R}^{n_v \times n_v}$ are supposed to be large and sparse. They represent the mass matrix and system matrix, respectively. $G \in \mathbb{R}^{n_v \times n_p}$ is a full rank matrix represents the discretized gradient and $B \in \R^{n_v \times n_b}$ is the input matrix. The system matrix $A \in \mathbb{R}^{n_v \times n_v}$ can be decomposed  as follows  $$A=-\dfrac{1}{\text{Re}}L-K-R.$$

\noindent More precisely, $-L\textbf{v}$ represents the discrete Laplacien $\Delta \, \textbf{v} $, $-K\textbf{v}$ is the discrete convection resulting from $(w \cdot \nabla) \textbf{v}$ and $-R \textbf{v}$ refers to the discrete reaction process of $(\textbf{v} \cdot \nabla) w$. In the Stokes system, $A$ is symmetric  negative definite matrix since there is no role to the matrices $K$ and $R$. A computational methods based on Krylov projection techniques and interpolatory projection to built a reduced system to a Stokes system have been established respectively in \cite{Ham,Guger}. We add to the system (\ref{DAE1}) an output function given by 
$$y(t)= C \bv(t),$$ 
where $y(t)$ is the output vector and $C^T\in \R^{n_v \times n_c}$ is the output matrix that measures velocity behaviour using information from internal nodes \cite{benner15}. The system (\ref{DAE1}) can be rewritten in a \textit{compact} form 
\begin{equation}
	\label{syscompact}
	\left\{\begin{array}{cclll}
		\underbrace{\left[\begin{array}{cccc}
				M &  0\\
				0  & 0\\
			\end{array}\right]}_{\bf M}\left[\begin{array}{ccc}
			\dot{\bv}(t)\\
			\dot{\textbf{p}}(t)\\
		\end{array}\right] & = &\underbrace{\left[\begin{array}{cccc}
				A &  G\\
				G^T  & 0\\
			\end{array}\right]}_{\bf A}\left[\begin{array}{ccc}
			\bv(t)\\
			\textbf{p}(t)\\
		\end{array}\right] +\left[\begin{array}{ccc}
			B\\
			0\\
		\end{array}\right]{\bf u}(t),\\[0.3cm]
		\textbf{y}(t)  &  = & \left[\begin{array}{ccc}
			C & 0\\
		\end{array}\right]\left[\begin{array}{ccc}
			\bv(t)\\
			\bp(t)\\
		\end{array}\right], \end{array}\right.
\end{equation}
and we call it a descriptor system since the matrix $\bf M$ is singular. It uses  the following matrix-pencil
\begin{equation}
	\label{pencil}
	\left(\left[\begin{array}{cccc}
		A &  G\\
		G^T  & 0\\
	\end{array}\right], \left[\begin{array}{cccc}
		M &  0\\
		0  & 0\\
	\end{array}\right]\right). 
\end{equation}
This matrix pencil has $n_v-n_p$ finite eigenvalues $\lambda_i \in \C \setminus {0}$ and $2n_p$ infinite eigenvalues $\lambda_{\infty} = \infty$, see \text{Theorem 2.1} in \cite{eigenvalues}. The system (\ref{syscompact}) is known as an index-2 descriptor dynamical system, for more details about the index of differential algebraic equation, see \cite{index}. For a Reynolds number ($Re \ge 300$), some eigenvalues of the matrix pencil ($\bf A, \bf M$) lie in $\C^+$, see \cite{benner15}.\\
Next, we present a whole process that allows us to reduce such systems. We describe a model reduction technique via a Krylov subspace-based 
method in order to construct an efficient reduced order system to (\ref{syscompact}) that has nearly the same response characteristics. To guarantee a well processing
of our suggested method, we need to establish a transformation of the system
(\ref{syscompact}) into an ordinary differential equations (ODEs).

\subsection{Deriving the ODE system}
\label{trans-ode}
We first eliminate the discrete pressure \textbf{p} from (\ref{eq1}) using  the following projection operator 
$$  \Pi = I_n - G\,(G^T  M^{-1}  G)^{-1} \, G^T M^{-1}   \in \mathbb{R}^{n_v\times n_v}.$$
It  is easy to check that
$$(\Pi^T)^2 = \Pi^T, \quad \Pi^2 = \Pi,\quad \Pi \, G = 0, \quad \Pi \, M = M \, \Pi^T \, {\rm and}\,  M^{-1} \, \Pi = \Pi^{T} M^{-1}. $$
The projection $\Pi^T$ is an $M$-orthogonal projection where 
for $x,\, y\in \R^{n_v}$ and $M\in \R^{n_v \times n_v}$,  the  $M$-inner product is defined by 
\begin{equation*}
	<x,y>_M = (x,My) = y^T Mx \quad (M \,\text{is a symmetric positive definite matrix}).
\end{equation*}
Notice that 
$$\text{null}(\Pi^T)=\text{range}(M^{-1} G) \quad and \quad \text{range}(\Pi^T)=\text{null}(G^T).$$
By using all these properties  we can show that 
\begin{equation}
	\label{prop1}
	0= G^T \textbf{v}(t) \qquad \text{if and only if } \qquad \textbf{v}(t)=\Pi^T \textbf{v}(t).
\end{equation}

\noindent Multiplying (\ref{eq1}) by $G^T\,M^{-1}$ and using (\ref{cond1}), the term $\bp$ can be expressed as follows 
$$\bp(t) = -(G^T\, M^{-1} \, G)^{-1}G^T\, M^{-1}\,A\bv(t)-(G^T\, M^{-1} \, G)^{-1} G^T\,M^{-1}\,B \textbf{u}(t).$$
Replacing $\textbf{p}$ in (\ref{eq1}) and multiplying by $\Pi$ yields to the following projected system
\begin{subequations}
	\label{pisys}
	\begin{align}	
		\cM \, \dfrac{d}{dt}\bv(t) &= \cA \bv(t)+\cB\textbf{u}(t),\\
		\textbf{y}(t) &= \cC \bv(t).
	\end{align}
\end{subequations}

\noindent Where $\cA= \Pi\, A\,\Pi^T, \, \cM=\Pi\, M\,\Pi^T, \, \cB=\Pi\, B$ and $\cC=C\,\Pi^T$. 
Since the matrix-pencil given by (\ref{pencil}) has $n_v-n_p$ finite eigenvalues \cite{eigenvalues}, a decomposition of $\Pi$ can be made by employing the thin singular value decomposition which leads to the following decomposition 
$$\Pi =\Theta_l \Theta_r^T,$$
where $\Theta_l, \, \Theta_r \in \R^{n_v\times (n_v-n_p)}$, are  full rank matrices satisfying
$$\Theta_l^T \Theta_r = \text{I}_{(n_v-n_p)}.$$
By inserting this decomposition into (\ref{pisys}) and considering a new variable $\tilde{\bv}(t)= \Theta_l^T \bv(t)$ with $\Theta_r\tilde{\bv}(t)= \Theta_r\Theta_l^T \bv(t)=\Pi^T \bv(t)=\bv(t)$, 
we get the following ODE system
\begin{subequations}
	\label{thetarsys}
	\begin{align}
		\label{mthetasys}
		M_{\Theta} \, \dfrac{d}{dt}\tilde{\bv}(t) &= A_{\Theta} \tilde{\bv}(t)+ B_{\Theta}\textbf{u}(t),\\
		y(t) &= C_{\Theta}\tilde{\bv}(t),
	\end{align}		
\end{subequations}
where $M_{\Theta}= \Theta_r^T M \Theta_r, \,A_{\Theta}= \Theta_r^T A \Theta_r \in \R^{(n_v-n_p) \times (n_v-n_p)}, \, B_{\Theta}= \Theta_r^T B \in \R^{(n_v-n_p) \times n_b},$ $C_{\Theta}= C \Theta_r \in \R^{n_c \times (n_v-n_p)}$. The matrix  $M_{\Theta}$ is non-singular due to the fact that $M$ is symmetric and positive definite. Notice  that the three systems (\ref{DAE1}), (\ref{pisys}) and (\ref{thetarsys}) are equivalent in the sense that their finite spectrum is the same \cite{spectrum} and also they realize the same transfer function. Before proving this result we give the definition of a transfer function associated  to the dynamical system (\ref{thetarsys}), to this end, we apply the Laplace transform given by
$$L(f)(s) :=\int_{0}^{\infty}e^{-st}f(t)dt,$$
to the system (\ref{thetarsys}), then we get the new system in the frequency domain
\begin{equation*}
	\left\{ \begin{array}{lll} sM_{\Theta}\, \widetilde{\bf V}(s) &=& A_{\Theta}\,\widetilde{\bf V}(s)+B_{\Theta}\,U(s), \\ Y(s) &=& C_{\Theta}\,\widetilde{\bf V}(s). \end{array} \right.
\end{equation*}
Where $\widetilde{\bf V}(s), \, {\bf U}(s)$ and ${\bf Y}(s)$ are the Laplace transform of $\tilde{\bv}(t), \, {\bf u}(t)$ and $y(t)$ respectively. By eliminating $\widetilde{\bf V}(s)$ from the two equations, we obtain
$${\bf Y}(s)=F_{\Theta}(s)\,{\bf U}(s),$$
where
\begin{equation}
	F_{\Theta}(s)= C_{\Theta} (sM_{\Theta}- A_{\Theta})^{-1} B_{\Theta},
\end{equation} 
is the transfer function associated to the system (\ref{thetarsys}).

\begin{remark}
	Let $F_m$ be the transfer function associated to the reduced system. In order to measure the accuracy of the resulting reduced system, we have to compute the error $\|F_{\Theta}-F_m\|$ with respect to a specific norm. This error can also be used to know how the response of the reduced system is close to that of the original one since $\Vert {\bf Y}(s)-{\bf Y}_m(s) \Vert \le \Vert F_{\Theta}(s)-F_m(s) \Vert \, \Vert {\bf U}(s) \Vert$.
\end{remark}

Denote by $X=\Theta_r (sM_{\Theta}- A_{\Theta})^{-1} B_{\Theta}$, then $F_{\Theta}= C X$. In addition, $X$ satisfies
$$B_{\Theta} =	(sM_{\Theta}- A_{\Theta}) \Theta_l^T X, $$	or equivalently, 
	$$\Pi B  =	\Pi (sM-A) \Pi^T X.$$
Due to the facts that $range(\Pi^T)= null(G^T)$ and $G$ is of full rank, we can verify that 

	\begin{align*}
		\begin{bmatrix}
			s M-A & -G \\ 
			-G^T & 0
		\end{bmatrix}\begin{bmatrix}
			X \\ 
			\star
		\end{bmatrix} = \begin{bmatrix}
			B \\ 
			0
		\end{bmatrix}. \\
	\end{align*}
	In fact, the relation $range(\Pi^T)= null(G^T)$ guarantees 
	$$G^TX=0,$$
	and the full rank $G$ leads to 
	$$\star = (G^TG)^{-1} G^T [(s M-A)X-B],$$
	thus, the desired result 
	\begin{align}
		\label{tf}
		F_{\Theta}(s) = C \, X &=\begin{bmatrix}
			C & 0
		\end{bmatrix}\begin{bmatrix}
			X\\
			\star
		\end{bmatrix}
		= \begin{bmatrix}
			C & 0
		\end{bmatrix}\begin{bmatrix}
			s M-A & -G \\ 
			-G^T & 0
		\end{bmatrix}^{-1}\begin{bmatrix}
			B \\ 
			0
		\end{bmatrix} = F(s),
\end{align}
where $F(s)$ is the transfer function associated to the original system (\ref{DAE1}). The technique used here  allows us to solve a saddle point problem instead of solving a linear system depending on the dense matrix $\Pi$ and its $\Theta$-decomposition as  established earlier  in \cite{bt}.

\begin{remark}
	We  notice that instead of reducing the original system (\ref{syscompact}), we can reduce the ODE system (\ref{thetarsys})  since it has  the same transfer functions as it is shown in (\ref{tf}). 
\end{remark}

The  matrices involved in  (\ref{thetarsys}) are  dense due to the  projector $\Pi$ and its decomposition and  that is why we need a strategy to avoid using direct computations with these matrices. In the next subsection, we show how to construct a reduced order system to (\ref{thetarsys}) by using the structure of the original system (\ref{syscompact}) without requiring any explicit computation of the  dense matrices ($M_{\Theta}, A_{\Theta}, B_{\Theta}, C_{\Theta}$), and this  leads to a considerable  saving of  cost and storage. Our calculations involve the implicit use of the system (\ref{thetarsys}) and this implies solving  saddle point problems. Details are given in the next section.

\section{A model reduction method to a descriptor index-2 dynamical system}\label{sec3}
Our goal is to find a reduced system to (\ref{thetarsys}) since it realizes the same transfer function of (\ref{syscompact}) as we mentioned before. This new system can be constructed using a projection technique onto an  extended block Krylov subspace that is defined in the following subsection.

\subsection{The extended block  Arnoldi algorithm}
Multiplying from the left of the first equation of the system (\ref{thetarsys}) by the inverse of $M_{\Theta}$ gives the following system which will be called the {\it $\Theta$-system}
\begin{subequations}
	\label{thetasys}
	\begin{align}
		\dfrac{d}{dt}\tilde{\bv}(t) &= M_{\Theta}^{-1}A_{\Theta} \tilde{\bv}(t)+ M_{\Theta}^{-1}B_{\Theta}\textbf{u}(t),\\
		y(t) &= C_{\Theta}\tilde{\bv}(t).
	\end{align}
\end{subequations}
The extended block Krylov subspace associated to the pair $(M_{\Theta}^{-1}A_{\Theta},M_{\Theta}^{-1}B_{\Theta})$ is defined  as follows 
\begin{align*}
	\bK^{ext}_m(M_{\Theta}^{-1}A_{\Theta},M_{\Theta}^{-1}B_{\Theta})&={\tt Range}([ (M_{\Theta}^{-1}A_{\Theta})^{-m}(M_{\Theta}^{-1}B_{\Theta}),\ldots,(M_{\Theta}^{-1}A_{\Theta})^{-1}(M_{\Theta}^{-1}B_{\Theta})\\  & (M_{\Theta}^{-1}B_{\Theta}),(M_{\Theta}^{-1}A_{\Theta}) (M_{\Theta}^{-1}B_{\Theta}), \ldots, (M_{\Theta}^{-1}A_{\Theta})^{m-1} (M_{\Theta}^{-1}B_{\Theta})]).
\end{align*}
The extended block Arnoldi algorithm for the pair  $(M_{\Theta}^{-1}A_{\Theta},M_{\Theta}^{-1}B_{\Theta})$ is summarized in the following algorithm.  

\begin{algorithm}[H]
	\caption{The extended block Arnoldi algorithm associated to the {\it $\Theta$ system}}\label{ext1}
	\begin{itemize}
		\item Inputs: $ M_{\Theta} \in \R^{(n_v-n_p) \times (n_v-n_p)}, \, \,A_{\Theta} \in \R^{(n_v-n_p) \times (n_v-n_p)}$, $ B_{\Theta} \in \R^{(n_v-n_p) \times n_b}$ and $m$.
		\item Compute  $[\cV_1^b, \Lambda]= {\tt qr}([ M_{\Theta}^{-1}B_{\Theta}, (M_{\Theta}^{-1}A_{\Theta})^{-1}\, M_{\Theta}^{-1}B_{\Theta}])$.
		\item For $j=1,\ldots,m $
		\begin{enumerate}
			\item Set $\cV_j^{(1)}$: first $n_b$ columns of $\cV_j^b$; $\cV_j^{(2)}$: second $n_b$ columns of $\cV_j^b$.
			\item  $\widetilde \cV_{j+1} =  [(M_{\Theta}^{-1}A_{\Theta})\,\cV_j^{(1)}, (M_{\Theta}^{-1}A_{\Theta})^{-1}\,\cV_j^{(2)}]$.
			\item  Orthogonalize $\widetilde \cV_{j+1}$ with respect to  $\cV_1^b,\ldots,\cV_j^b$ to get $\cV_{j+1}^b$, i.e.,\\
			\hspace*{0.5cm} for $ i=1,2,\ldots,j $ \\
			\hspace*{1cm} $ H_{i,j} =  (\cV_i^b)^{T} \,\widetilde \cV_{j+1} $.\\
			\hspace*{1cm} $ \widetilde \cV_{j+1} = \widetilde \cV_{j+1} - \cV_i^b\,H_{i,j}  $.\\
			\hspace*{0.5cm} end for
			\item   $ [\cV_{j+1}^b, \; H_{j+1,j}] = QR(\widetilde \cV_{j+1})$.
			\item $\cV_{j+1} = [\cV_j, \; \cV_{j+1}^b]$.
		\end{enumerate}
		End For.
	\end{itemize}${}$
	\label{thetaalgo}
\end{algorithm} 
The extended block Arnoldi algorithm allows us to construct an orthonormal basis of  $\bK^{ext}_m(M_{\Theta}^{-1}A_{\Theta},M_{\Theta}^{-1}B_{\Theta})$ formed by the columns of $\{\cV_1^b,\ldots,\cV_m^b\}$, where $\cV^b_j$ for ($j=1, \ldots, m$) are $(n_v-n_p) \times 2n_b$ matrices. We also have some classical algebraic properties given in the following proposition.

\begin{proposition}
	Let $\cV_m =[\cV_1^b,\ldots, \cV_m^b] \in \R^{2mn_b \times 2mn_b}$ be the matrix generated using the extended block Arnoldi Algorithm \ref{thetaalgo}  to the pairs $(M_{\Theta}^{-1}A_{\Theta},M_{\Theta}^{-1}B_{\Theta})$, $\bT_m= \cV_m^T \, M_{\Theta}^{-1}A_{\Theta} \, \cV_m$. Then we have the following results
	\begin{align}
		\label{eq3.1}
		M_{\Theta}^{-1}  A_{\Theta} \, \cV_m &= \cV_{m+1} \, \overline \bT_m  \\ 
		&=  \cV_m \, \bT_m + \cV_{m+1}^b \, T_{m+1,m} E_m^T,
	\end{align}
	where $T_{m+1,m}$ is the last $2n_b \times 2n_b$ block of $\overline{\bT}_m \in \R^{2(m+1)n_b\times 2mn_b}$ and $E_m^T$ is the last $2n_b$ columns of the identity matrix $I_{2mn_b}$.
\end{proposition} 
\begin{proof}
	Using the fact that
	 $$M_{\Theta}^{-1}  A_{\Theta} \bK^{ext}_m(M_{\Theta}^{-1}A_{\Theta},M_{\Theta}^{-1}B_{\Theta}) \subset \bK^{ext}_{m+1}(M_{\Theta}^{-1}A_{\Theta},M_{\Theta}^{-1}B_{\Theta}), $$
	 and the orthogonality of $\cV_m$,  there exists a matrix $L$ such that 
	\begin{align}
		\label{eqalgtheta}
		M_{\Theta}^{-1}  A_{\Theta} \, \cV_m = \cV_{m+1} \, L.
	\end{align}
	It has been shown that $\bT_m$ is an upper block Hessenberg matrix in \cite{heyouni,simoncini} and also that $\bT_m$ can be computed directly from the columns of the upper block Hessenberg matrix $\bH_m$ generated by Algorithm \ref{newextd}. Since $\cV_{m+1}=[\cV_m, \cV_{m+1}^b],$ we have 
	\begin{eqnarray*}
		\bT_{m+1} &=& \cV_{m+1}^T \, M_{\Theta}^{-1}A_{\Theta} \, \cV_{m+1} \\
		&=& \begin{bmatrix}
			\cV_m^T \, M_{\Theta}^{-1}A_{\Theta} \, \cV_m &  \cV_m^TM_{\Theta}^{-1}A_{\Theta} \, \cV_{m+1}^b \\ 
			(\cV_{m+1}^b)^T \, M_{\Theta}^{-1}A_{\Theta} \, \cV_m & (\cV_{m+1}^b)^T M_{\Theta}^{-1}A_{\Theta}\cV_{m+1}^b 
		\end{bmatrix} \\
		&=& \begin{bmatrix}
			\bT_m & \cV_m^TM_{\Theta}^{-1}A_{\Theta} \, \cV_{m+1}^b \\ 
			(\cV_{m+1}^b)^T \, M_{\Theta}^{-1}A_{\Theta} \, \cV_m & (\cV_{m+1}^b)^T M_{\Theta}^{-1}A_{\Theta}\cV_{m+1}^b 
		\end{bmatrix}.
	\end{eqnarray*}
	We know that $\bT_{m+1}$ is also un upper block Hessenberg matrix, then $$T_{m+1,m}\, E_m^T = (\cV_{m+1}^b)^T \, M_{\Theta}^{-1}A_{\Theta} \, \cV_m,$$ and 
	$$\overline \bT_m  = \cV_{m+1}^T M_{\Theta}^{-1}A_{\Theta} \cV_m =\begin{bmatrix}
		\bT_m\\ 
		T_{m+1,m} E_m^T
	\end{bmatrix} \in \R^{2(m+1)n_b \times 2mn_b}.$$
	Multiplying by  $\cV_{m+1}^T$ from the left of (\ref{eqalgtheta}), we obtain $\overline \bT_{m} = L$. As a consequence we get the desired result
	\begin{align*}
		M_{\Theta}^{-1}  A_{\Theta} \, \cV_m &= \cV_{m+1} \, \overline \bT_m \label{eq11} \\ 
		&=[\cV_m, \cV_{m+1}^b]\begin{bmatrix}
			\bT_m\\ 
			T_{m+1,m} E_m^T
		\end{bmatrix} \\
		&=  \cV_m \, \bT_m + \cV_{m+1}^b \, T_{m+1,m} E_m^T.
	\end{align*}
\end{proof} 

After constructing the matrix $\cV_m$ corresponding to the basis of the extended block Krylov subspace  $\bK^{ext}_m(M_{\Theta}^{-1}A_{\Theta},M_{\Theta}^{-1}B_{\Theta}) $, 
we can now built the reduced system by considering the approximation $\tilde{\bv}(t) \approx \cV_m \bv_m(t)$ and by replacing in  (\ref{thetarsys}), and then imposing  the Petrov-Galerking condition,  we obtain the following projected reduced order dynamical system 
\begin{equation}
	\left\{ \begin{array}{lll} \cM_m\dot{\bv}_m(t) &=& \cA_m\,\bv_m(t)+\cB_m\,{\bf u}(t), \\ y_m(t) &=& \cC_m\,\bv_m(t), \end{array} \right.
\end{equation}
with the associated transfer function $F_m(s)=\cC_m(s\cM_m-\cA_m)^{-1}\cB_m$, where $\cM_m=\cV_m^T \cM \cV_m, \, \, \cA_m=\cV_m^T \cA \cV_m \in \R^{2mn_b \times 2mn_b}$ and  $\cB_m=\cV_m^T \cB \in \R^{2mn_b \times n_b}, \,\, \cC_m= C \cV_m \in \R^{n_c \times 2mn_b }$. \\
As we mentioned before, the explicit computation of $\cV_m$ is prohibitive in our approach since the $j$-th block $\cV_j^b$ of $\cV_m$ relies on $\Theta_r$, which will make our calculations infeasible due to the density of the $\Theta$-decomposition of the projection $\Pi$. In what follows, we describe an appropriate process to get a reduced system to (\ref{thetasys}) by avoiding an explicit computation of $\cV_m$.\\
The main computational issue when we apply the extended block Arnoldi Algorithm \ref{thetaalgo} to the pair $(M_{\Theta}^{-1}A_{\Theta},M_{\Theta}^{-1}B_{\Theta})$  is to compute   blocks of the form  
	\begin{align}
		\label{block1}
		\widetilde{\cV}_1 &= [M_{\Theta}^{-1}B_{\Theta}, (M_{\Theta}^{-1}A_{\Theta})^{-1} M_{\Theta}^{-1}B_{\Theta}]\\
		&=[\widetilde{\cV}_1^{(1)},\widetilde{\cV}_1^{(2)}],
	\end{align}
	and for $j=1,\ldots,m$,
	\begin{align}
		\label{block2}
		\widetilde{\cV}_{j+1} &= [(M_{\Theta}^{-1}A_{\Theta})\,\cV_j^{(1)}, (M_{\Theta}^{-1}A_{\Theta})^{-1}\,\cV_j^{(2)}]\\
		&=[\widetilde{\cV}_{j+1}^{(1)},\widetilde{\cV}_{j+1}^{(2)}],
	\end{align}
	where $\cV_j^{(1)}$ and $\cV_j^{(2)}$ are the first and second $n_b$ columns of $\cV_j^b$,  respectively. Our strategy consists in  reformulating those blocks onto new ones  without an explicit calculation of $\Theta_r$. We set $\widetilde{\bV}_m = \Theta_r \widetilde{\cV}_m \in \R^{n_v \times 2mn_b}$ where $\widetilde{\cV}_m =[\widetilde{\cV}_1,\ldots, \widetilde{\cV}_m] \in \R^{(n_v-n_p) \times 2mn_b}$ and $\widetilde{\bV}_m =[\widetilde{V}_1,\ldots, \widetilde{V}_m] \in \R^{n_v\times 2mn_b}$ satisfying  
	\begin{align}
		\label{Pitvm_old}
		\Pi^T \widetilde{\bV}_m= \Theta_r \Theta_l^T \widetilde{\bV}_m=\Theta_r \widetilde{\cV}_m = \widetilde{\bV}_m,
	\end{align}
	We set again $\bV_m = \Theta_r \cV_m \in \R^{n_v \times 2mn_b}$. All the $j$-th block $V_j \in \R^{n_v \times 2n_b}$ of $\bV_m$ are computed in an appropriate way, which means that we do not include the matrix $\Theta_r$ in our computation and also not the block $\cV_j^b$. Details are given in Algorithm \ref{newextd}. 
	\begin{align}
		\label{Pitvm}
		\Pi^T \bV_m= \Theta_r \Theta_l^T \bV_m=\Theta_r \cV_m = \bV_m,
	\end{align}
	The result (\ref{Pitvm_old}) confirms that $G^T \widetilde{\bV}_m=0$ as it is shown in (\ref{prop1}), and  consequently we obtain the following relations
	\begin{itemize}
		\item 	$M_{\Theta}^{-1}B_{\Theta} = \widetilde{\cV}_1^{(1)}$,    
		\item 	$M_{\Theta} \widetilde{\cV}_1^{(1)}=B_{\Theta}$, 
		\item 	$\Theta_r^T M \Theta_r \widetilde{\cV}_1^{(1)} =\Theta_r B$, 
		\item 	$\Pi M \Pi^T \widetilde{V}_1^{(1)} =\Pi B$, 
		\item 	$\Pi (M  \widetilde{V}_1^{(1)}-B) =0$, 
		\item 	$(M  \widetilde{V}_1^{(1)}-B) \in null(\Pi)=range(G)$.
	\end{itemize}
	Then,  the first $n_b$ block-column  $\widetilde{V}_{1}^{(1)}$ of $\widetilde{V}_{1} \in \R^{n_v \times 2n_b}$ can be computed by solving the following saddle point problem
	$$\begin{bmatrix}
		M & G \\ 
		G^T & 0
	\end{bmatrix} \begin{bmatrix}
		\widetilde{V}_1^{(1)} \\ 
		\star
	\end{bmatrix} = \begin{bmatrix}
		B \\ 
		0
	\end{bmatrix}.$$
	The same process can be used  to get $\widetilde{V}_1^{(2)}$ by starting from the following linear system $$(M_{\Theta}^{-1}A_{\Theta})^{-1}\, M_{\Theta}^{-1}B_{\Theta} =\widetilde{\cV}_1^{(2)}.$$ After that, one can use the {\tt qr} function (in MATLAB) to find the block $V_1=[V_1^{(1)}, V_1^{(2)}] \in \R^{n_v \times 2n_b}$ as described in Algorithm \ref{newextd}.
	To get the first $n_b$ block-column $\widetilde{V}_{j+1}^{(1)}$ of $\widetilde{V}_{j+1}$, we  use the following steps 
	\begin{itemize}
		\item $(M_{\Theta}^{-1}A_{\Theta})\cV_{j}^{(1)} = \widetilde{\cV}_{j+1}^{(1)},$
		\item $M_{\Theta} \widetilde{\cV}_{j+1}^{(1)}=A_{\Theta}\cV_{j}^{(1)},$ 
		\item $\Theta_r^T M \Theta_r \widetilde{\cV}_{j+1}^{(1)} =\Theta_r^T A \Theta_r \cV_{j}^{(1)},$ 
		\item $\Pi M \Pi^T \widetilde{V}_{j+1}^{(1)} =\Pi A V_j^{(1)},$ 
		\item $\Pi (M  \widetilde{V}_{j+1}^{(1)}-A V_j^{(1)}) =0,$  
		\item $(M  \widetilde{V}_{j+1}^{(1)}-A V_j^{(1)}) \in null(\Pi)=range(G).$
	\end{itemize}
	\noindent Then we have to solve the following saddle point problem
	$$\begin{bmatrix}
		M & G \\ 
		G^T & 0
	\end{bmatrix} \begin{bmatrix}
		\widetilde{V}_{j+1}^{(1)} \\ 
		\star
	\end{bmatrix} = \begin{bmatrix}
		A V_j^{(1)} \\ 
		0
	\end{bmatrix}.$$
In the same manner, we can compute the last $n_b$ column $\widetilde{V}_{j+1}^{(2)}$ of $\widetilde{V}_{j+1}$ by starting from this linear system $(M_{\Theta}^{-1}A_{\Theta})^{-1}\,\cV_j^{(2)} =\widetilde{\cV}_{j+1}^{(2)} $ and following the same previous process.\\ 
After showing how to compute the block vectors (\ref{block1}) and (\ref{block2}) without computing neither the matrix $\cV_m$ corresponding to the orthonormal basis of $\bK^{ext}_m(M_{\Theta}^{-1}A_{\Theta},M_{\Theta}^{-1}B_{\Theta})$ nor $\Theta$-decomposition of $\Pi$, we can now present the new extended block Arnoldi algorithm based only on the sparse system matrices of the index-2 system. Here, we have to mention that this algorithm is based on a Gram-Shmidt orthogonalization process, which reconstructs the blocks $\{V_1,\ldots, V_m\}$, such that their columns form an orthonormal matrix $\bV_m$ as described in Algorithm \ref{newextd} step 3.c. This matrix will be used in order to get an efficient reduced system to the index-2 original one (\ref{DAE1}). Details are given in the next subsections. We summarize  all these steps in the following algorithm. 
\begin{algorithm}[H]
	\caption{The extended block Arnoldi algorithm associated to the index-2 system}
	\begin{itemize}
		\item[] Inputs: $ M\in \R^{n_v \times n_v}, \,A\in \R^{n_v \times n_v}, \, G\in \R^{n_v \times n_p},\,B\in \R^{n_v \times n_b}$ and $m$. 
		\item[1.] solving the first saddle point problems 
		$$\begin{bmatrix}
			M & G \\ 
			G^T & 0
		\end{bmatrix} \begin{bmatrix}
			\widetilde{V}_1^{(1)} \\ 
			\star
		\end{bmatrix} = \begin{bmatrix}
			B \\ 
			0
		\end{bmatrix}, \begin{bmatrix}
			A & G \\ 
			G^T & 0
		\end{bmatrix} \begin{bmatrix}
			\widetilde{V}_1^{(2)} \\ 
			\star
		\end{bmatrix} = \begin{bmatrix}
			B \\ 
			0
		\end{bmatrix},$$
		\item[2.] Compute  $[V_1, \Lambda]= {\tt qr}\left( [\widetilde{V}_1^{(1)},\widetilde{V}_1^{(2)}]\right)$, $ \bV_1 = [V_1].$
		\item [3.] For $j=1,\ldots,m $
		\begin{enumerate}
			\item[a.] Set $V_j^{(1)}$: first $n_b$ columns of $V_j$; $V_j^{(2)}$: second $n_b$ columns of $V_j$.
			\item[b.]  $\widehat V_{j+1} = \left(\begin{bmatrix}
				M & G \\ 
				G^T & 0
			\end{bmatrix} \begin{bmatrix}
				\widetilde{V}_{j+1}^{(1)} \\ 
				\star
			\end{bmatrix} = \begin{bmatrix}
				A\,V_j^{(1)} \\ 
				0
			\end{bmatrix}, \begin{bmatrix}
				A & G \\ 
				G^T & 0
			\end{bmatrix} \begin{bmatrix}
				\widetilde{V}_{j+1}^{(2)} \\ 
				\star
			\end{bmatrix} = \begin{bmatrix}
				M\,V_j^{(2)} \\ 
				0
			\end{bmatrix}\right)$.
			\item[c.]  Orthogonalize $\widehat V_{j+1}$ with respect to  $V_1,\ldots,V_j$ to get $V_{j+1}$, i.e.,\\
			\hspace*{0.5cm} for $ i=1,2,\ldots,j $ \\
			\hspace*{1cm} $ H_{i,j} =  (V_i)^{T} \,\widehat V_{j+1} $; \\
			\hspace*{1cm} $ \widehat V_{j+1} = \widehat V_{j+1} - V_i\,H_{i,j}  $; \\
			\hspace*{0.5cm} end for
			\item[d.]   $ [V_{j+1}, \; H_{j+1,j}] = QR(\widehat V_{j+1})$.
			\item[e.] $\bV_{j+1} = [\bV_j, \; V_{j+1}]$.
		\end{enumerate}
		End For.
	\end{itemize}${}$
	\label{newextd}
\end{algorithm} 
\noindent As we noticed, the main steps of Algorithm \ref{newextd} is the solution of a saddle-point problems of $n_v+n_p$ dimension in Step 1 and in Step 3.b, and we are interesting only in the first $n_v$ rows. At each iteration, a direct solver $"\backslash", \, \text{ a built-in function on MATLAB}$, is used to solve these saddle point problems. The new vector $V_{j+1}$ of the matrix $\bV_m$ can be computed via the  Gram-Shmidt process as we explain in the Step 3.c. The $"\star"$ refers to an $n_p \times n_b$ block  that is not taken into account. 
After $m$ steps of Algorithm \ref{newextd}, we get an orthonormal matrix  $\bV_m = \left [V_1,V_2,\ldots,V_m \right ] \in \R^{n_v \times 2mn_b}$ with $V_i \in \R^{n_v \times 2n_b}$. This algorithm built also  an upper block Hessenberg matrix $\bH_m \in \R^{2mn_b \times 2mn_b}$ whose non zero blocks are the $H_{i,j}$. Notice  that each submatrix $H_{i,j}$ ($1 \le i \le j \le m $) is of order $2n_b \times 2n_b$. 
A similar algebraic relations to the one given by (\ref{eq3.1}) can be derived  using only the sparse matrices $M,A$ and also the matrix $\bV_m$ generated by Algorithm \ref{newextd}. We present this result in the following proposition. 
\begin{proposition}
	\label{Tm}
	Let $\bV_m \in \R^{n_v \times 2mn_b}$ and $\overline \bT_m \in \R^{2(m+1)n_b \times 2mn_b}$   be the  orthonormal matrix and the upper block Hessenberg matrix generated by Algorithm\ref{newextd}, respectively. Then we have 
	\begin{align*}
		M^{-1} \Pi\, A \bV_m &=\bV_{m+1} \, \overline \bT_m \\
		&= \bV_m \bT_m + V_{m+1} T_{m+1,m} E_m^T,
	\end{align*}
	where 	$\Pi$ is the projection matrix defined earlier.
\end{proposition}
\begin{proof}
	Multiplying  from the left the relation  (\ref{eq3.1}) by  $\Theta_r$, and using the fact that $\bV_m = \Theta_r \cV_m$, we get	
	\begin{align}
		\Theta_r M_{\Theta}^{-1} \Theta_r^T A\Theta_r \, \cV_m &= \Theta_r\cV_{m+1} \, \overline \bT_m, \\ 
		\label{eq12}\Theta_r M_{\Theta}^{-1}\Theta_r^T A\bV_m &=\bV_{m+1} \, \overline \bT_m.
	\end{align}
	On the other hand, we know that $\Pi M= M \Pi^T$ by definition of $\Pi$, and by using the fact that  $\Pi M \Theta_r= M\Theta_r$ by  the $\Theta$-decomposition, we obtain the following relations
	\begin{align*}
		\Pi M \Theta_r&= M\Theta_r,\\
		\Theta_l \Theta_r^T M \Theta_r &= M \Theta_r, \\
		\Theta_l M_{\Theta} &= M \Theta_r,\\
		M^{-1} \Theta_l &= \Theta_r M_{\Theta}^{-1}, \\
		\Theta_r M_{\Theta}^{-1}\Theta_r^T &= M^{-1} \Pi.
	\end{align*}	
	\noindent Replacing the last relation   in the formula (\ref{eq12}), we get the desired result.
\end{proof}

\noindent Notice that from Step 1 of Algorithm \ref{newextd}, we have
\begin{align}
	\label{step1algo}
	[V_1, \Lambda]= {\tt qr}\left([v,w]\right),
\end{align}
where $\Lambda \in \R^{2n_b \times 2n_b}$ is an upper triangular matrix defined by 
$$\Lambda = \begin{bmatrix}
	\Lambda^{(1,1)}	& \Lambda^{(1,2)} \\
	0	& \Lambda^{(2,2)}
\end{bmatrix}, 
$$
and $v,\, w$ are the solutions of the following saddle point problems
$$\begin{bmatrix}
	M & G \\ 
	G^T & 0
\end{bmatrix} \begin{bmatrix}
	v \\ 
	\star
\end{bmatrix} = \begin{bmatrix}
	B \\ 
	0
\end{bmatrix} \; and \; 
\begin{bmatrix}
	A & G \\ 
	G^T & 0
\end{bmatrix} \begin{bmatrix}
	w \\ 
	\star
\end{bmatrix} = \begin{bmatrix}
	B \\ 
	0
\end{bmatrix}.$$
We notice  that $$\begin{bmatrix}
	M & G \\ 
	G^T & 0
\end{bmatrix} \begin{bmatrix}
	v \\ 
	\star
\end{bmatrix} = \begin{bmatrix}
	B \\ 
	0
\end{bmatrix} \Leftrightarrow \Pi M \Pi^T v = \Pi B, \, \, (\text{with}\; \Pi^T v=v),$$ and from (\ref{step1algo}) we get  $$	[v,w]=[V_1^{(1)},V_1^{(2)} ] \, \begin{bmatrix}
	\Lambda^{(1,1)}	& \Lambda^{(1,2)} \\
	0	& \Lambda^{(2,2)}
\end{bmatrix},$$
thus
$$v=V_1^{(1)} \Lambda^{(1,1)}, $$
and then 
\begin{align}
	\label{b}
	\bV_m^T \, M^{-1} \, \Pi B = \bV_m^T V_1^{(1)} \Lambda^{(1,1)} = \begin{bmatrix} I_{n_b}  \\ 0_{n_b} \\ \vdots \\ 0_{n_b}  \end{bmatrix} \Lambda^{(1,1)}.
\end{align}
We have mentioned before that in order to reduce the original system (\ref{syscompact}), we can construct a reduced system from the  $\Theta$ system (\ref{thetarsys}) since they realize the same transfer function as it is shown in (\ref{tf}). At the iteration $m$, we approximate  $\tilde{\bv}(t)$ by $\cV_m \hat{\bf v}_m(t)$  where $\cV_m$ is the  matrix corresponding to the orthonormal basis of $\bK^{ext}_m(M_{\Theta}^{-1}A_{\Theta},M_{\Theta}^{-1}B_{\Theta})$. By injecting the approximation of $\tilde{\bv}(t)$ in the system (\ref{thetarsys}) and enforcing the Petrov-Galerkin condition, we get the following reduced system 
\begin{equation}
	\left\{ \begin{array}{lll} \dot{\hat{{\bf v}}}_m(t) &=& \cV_m^T \, M_{\Theta}^{-1}A_{\Theta} \, \cV_m\,\hat{\bf v}_m(t)+\cV_m^T M_{\Theta}^{-1}B_{\Theta}{\bf u}(t), \\ y_m(t) &=& C_{\Theta}\cV_m\,\hat{{\bf v}}_m(t). \end{array} \right.
\end{equation}
We know that $\bT_m= \cV_m^T \, M_{\Theta}^{-1}A_{\Theta} \, \cV_m $ which can be computed only from the upper block Hessenberg matrix $\bH_m$ generated by Algorithm \ref{newextd} as we mentioned before, also $C_{\Theta}\cV_m = C \Theta_r \cV_m= C \bV_m$, and by using the fact that $M_{\Theta}^{-1}B_{\Theta} \in \bK^{ext}_m(M_{\Theta}^{-1}A_{\Theta},M_{\Theta}^{-1}B_{\Theta})$ which confirms that $\cV_m\cV_m^T M_{\Theta}^{-1}B_{\Theta}=M_{\Theta}^{-1}B_{\Theta}$, then we can prove 
$$\cV_m^T M_{\Theta}^{-1}B_{\Theta} =\bV_m^T \, M^{-1} \, \Pi B.$$
Finally, we get the following reduced order LTI dynamical system 
\begin{equation}
	\label{sysredu1}
	\left\{ \begin{array}{lll} \dot{\hat{\bf v}}_m(t) &=& \bT_{m}\,\hat{\bf v}_m(t)+\bB_m{\bf u}(t), \\ y_m(t) &=& \bC_m\,\hat{\bf v}_m(t), \end{array} \right.
\end{equation}
where $\bB_m= \begin{bmatrix} I_{n_b}, \, 0_{n_b},\ldots,0_{n_b}  \end{bmatrix}^T \Lambda^{(1,1)} \in \R^{2mn_b \times n_b}$ as it is mentioned in (\ref{b}), and $\bC_m= C\bV_m \in \R^{n_c \times 2mn_b}$.\\
The reduced transfer function is given by 
$$F_m(s) = \bC_m(sI_{2mn_b}-\bT_m)^{-1} \bB_m.$$ 
Another way to construct a reduced system is by considering the system  (\ref{thetarsys}) without inverting the matrix $M_{\Theta}$. We again approximate $\tilde{\bv}(t)$ by $\cV_m \hat{{\bf v}}_m(t)$ where $\cV_m$ is a matrix described in the previous sections, then we get the following system 
\begin{equation}
	\left\{ \begin{array}{lll} \cV_m^T \, M_{\Theta}\cV_m\dot{\hat{\bf v}}_m(t) &=& \cV_m^T \, A_{\Theta} \, \cV_m\,\hat{\bf v}_m(t)+\cV_m^T B_{\Theta}{\bf u}(t), \\ y_m(t) &=& C_{\Theta}\cV_m\,\hat{{\bf v}}_m(t), \end{array} \right.
\end{equation}
using the fact that $\bV_m= \Theta_r \cV_m$, we get the following reduced system 
\begin{equation}
	\label{sysredu2}
	\left\{ \begin{array}{lll} \bM_m\dot{\hat{\bf v}}_m(t) &=& \bA_m\,\hat{{\bf v}}_m(t)+\bB_m{\bf u}(t), \\ y_m(t) &=& \bC_m\,\hat{{\bf v}}_m(t), \end{array} \right.
\end{equation}
with the associated reduced transfer function $$F_m(s) = \bC_m(s\bM_m-\bA_m)^{-1} \bB_m,$$
where $\bM_m= \bV_m^T M \bV_m,\, \bA_m= \bV_m^T A\bV_m, \, \bB_m= \bV_m^TB$ and $\bC_m=C\bV_m.$\\ In Algorithm \ref{newextd} we gave a description of the process to get the matrix $\bV_m$ without any explicit computation of $\cV_m$  or  $\Theta_r$. \\
\begin{remark}
	The two reduced dynamical  systems (\ref{sysredu1}) and (\ref{sysredu2}) are considered efficient reduced systems compared to the original one represented by the index-2 system (\ref{syscompact}), but numerically the first reduced system is more economical since its system matrices ($\bT_{m}\, \bB_m,\, \bC_m$) could be computed appropriately and without requiring matrix-vector products with $A$ and $M$ which is the case for the second reduced system represented by the system matrices ($\bM_m,\bA_m,\bB_m, \, \bC_m$).
\end{remark}

\section{Solving the LQR problem based on a Riccati feedback approach}
\label{sec4}
The linear quadratic regulator is a well-known classical method for constructing controlled feedback gains. This feedback allows the design of stable and efficient closed-loop systems. We used the transformation explained in Subsection \ref{trans-ode} that allows us to deal with an LQR problem governed by an ODE instead of an LQR problem governed by an DAE. Following that, a classical LQR theory can be applied to solve the new problem based on a Riccati feedback approach. The main issue with this approach is the solution of a generalized algebraic Riccati equation (GARe). We mentioned earlier that all calculations are performed using the structure of the original DAE system and not that of the ODE one due to the density of projection $\Pi$ and its $\Theta$-decomposition which can make our calculations infeasible.
\subsection{The LQR problem associated to the ODE system}
The LQR problem consists in  minimizing the following cost functional 
\begin{align} \label{cost}
	\cJ(\tilde{\bv},\textbf{u}(t)) :=\frac{1}{2}\int_{0}^{\infty} (\tilde{\bv}^TC_{\theta}^TC_{\theta}\tilde{\bv} + \textbf{u}(t)^T\textbf{u}(t)) \, \, \text{dt},
\end{align}
subject to the ODE system (\ref{thetarsys}) constraints defined earlier in Section \ref{sec2}. According to the LQR approach, the optimal control that minimizes the functional coast (\ref{cost}) subject to the  dynamical constraints (\ref{thetarsys}) is given  by 
\begin{align}
	\textbf{u}_{\star}(t) = -\underbrace{B_{\theta}^TX_{\theta} M_{\theta}}_{:=K_{\Theta}}\tilde{\bv},
\end{align} 
where $X_{\theta} \in \R^{(n_v-n_p) \times (n_v-n_p)}$ is the unique symmetric semi-definite positive stabilizing solution of the following generalized algebraic Riccati equation (GARe)
\begin{align}\label{riccati}
	\mR(X_{\Theta}) :=A_{\Theta}^TX_{\Theta}\,M_{\Theta} + M_{\Theta}\,X_{\Theta}\,A_{\Theta} -M_{\Theta}\,X_{\Theta}\,B_{\Theta}B_{\Theta}^TX_{\Theta}\,M_{\Theta}+ C_{\Theta}^TC_{\Theta} = 0.
\end{align}
The unique solution $X_{\Theta}$ can be computed using an extended block Krylov subspace method. This solution is the main ingredient to construct the feedback matrix $K_{\Theta} \in \R^{n_b \times (n_v-n_p)}$ that asymptotically stabilizes the ODE system (\ref{thetarsys}). However, solving the GARe (\ref{riccati}) is not recommended in our process due to the presence of $\Theta$-decomposition of the  projection $\Pi$. In the next subsections, we describe how to solve such an algebraic equation (\ref{riccati}) without using  the  $\Theta$-decomposition in our computations.
\subsection{Solving the generalized algebraic Riccati equation}

Our goal is to solve the GARe (\ref{riccati}) without any explicit computation of  the dense matrices ($M_{\Theta}, A_{\Theta}, B_{\Theta}, C_{\Theta}$). This statement intended to the fact that those matrices rely on $\Theta_r$, and a direct use of them can make our calculations impractical due to the density of $\Theta$-decomposition of the projection $\Pi$. 
A multiplication from the left and right of (\ref{riccati}) by $\Theta_l$ and $\Theta_l^T$,  respectively, gives the following result 
\begin{align*}
\Pi A^T \Theta_rX_{\Theta}\Theta_r^TM \Pi^T+ \Pi M \Theta_rX_{\Theta}\Theta_r^T A \Pi^T -\Pi M \Theta_rX_{\Theta}\Theta_r^TB B^T\theta_rX_{\theta}\theta_r^T &M\Pi^T\\
&+ \Pi C^TC \Pi^T = 0.
\end{align*}
Setting $X=\Theta_rX_{\Theta}\Theta_r^T$, and using  the fact that $\Pi M=M \Pi^T$, we get
$$\Pi A^T X \Pi M +  M \Pi^TX A \Pi^T - M\Pi^T XB B^TX \Pi M + \Pi C^TC \Pi^T = 0.$$
Since $X \Pi = \Theta_r X_{\Theta}\Theta_r^T \Theta_l \Theta_r^T= \Theta_rX_{\Theta}\Theta_r^T =X$, same to $\Pi^T X=X$, then we obtain the following final result 
\begin{align}\label{riccatinew}
	\Pi A^T X M + M X A \Pi^T - M XB B^TX  M + \Pi C^TC \Pi^T = 0.
\end{align}
If we  set $K=B^T X M$ the feedback matrix associated to (\ref{riccatinew}), then the relation between $K_{\Theta}$ and $K$ is given as $$K_{\Theta} = B^T\Theta_rX_{\Theta}\Theta_r^T M \Theta_r =B^T X M \Theta_r = K \Theta_r.$$

In what follows, we describe an appropriate process to compute  the unique solution $X=X^T \succeq 0$, by avoiding an explicit computation of $\Theta_r$ or the solution $X_{\Theta}$. Multiplying  GARe (\ref{riccati})  from the left and the  right by the inverse of $M_{\Theta}$, we get 
$$M_{\Theta}^{-1}A_{\Theta}^TX_{\Theta} + X_{\Theta}\,A_{\Theta} M_{\Theta}^{-1} -X_{\Theta}\,B_{\Theta}B_{\Theta}^TX_{\Theta}+ M_{\Theta}^{-1}C_{\Theta}^TC_{\Theta}M_{\Theta}^{-1} = 0.$$
Then we apply the extended block Arnoldi Algorithm \ref{thetaalgo} to the pair $(M_{\Theta}^{-1}A_{\Theta}^T,M_{\Theta}^{-1}C_{\Theta}^T).$ 
The same process described in Section \ref{sec3} is followed here. We set again $\bV_m = \Theta_r \cV_m \in \R^{n_v \times 2mn_c}$ satisfying 
\begin{align}
	\label{Pitvm2}
	\Pi^T \bV_m= \Theta_r \Theta_l^T \bV_m=\Theta_r \cV_m = \bV_m.
\end{align}
As we mentioned before,   the orthonormal matrix $\bV_m$ can be constructed using the Algorithm \ref{newextd} without any explicit computation of $\cV_m$ or $\Theta_r$. After $m$ iterations of the process,  we can use Proposition \ref{Tm} to prove that 
\begin{subequations}
	\label{equaalgTm}
	\begin{align}
		M^{-1} \Pi\, A^T \bV_m &=\bV_{m+1} \, \overline \bT_m \\
		&= \bV_m \bT_m + V_{m+1} T_{m+1,m} E_m^T.
	\end{align}
\end{subequations}
We seek for a low rank approximate solution to the GARe (\ref{riccatinew})
under the following form 
\begin{equation}
	\label{Xm}
	X_m=\bV_m Y_m \bV_m^T,
\end{equation} 
where $Y_m \in \R^{2mn_c \times 2mn_c}$ is the unique solution of a low-dimensional Riccati equation defined below. Replacing the approximation (\ref{Xm}) in the equation (\ref{riccatinew})  and multiplying from the left and right by the inverse of $M$, we obtain
\begin{align*}
	M^{-1}\Pi A^T \bV_m Y_m \bV_m^T +  \bV_m Y_m \bV_m^T A \Pi^T M^{-1} - \bV_m Y_m \bV_m^T&B B^T\bV_m Y_m \bV_m^T   \\
	&+ M^{-1}\Pi C^TC \Pi^TM^{-1} = 0,
\end{align*}
which gives 
$$\bT_{m} Y_m  +   Y_m \bT_{m}^T -  Y_m \bV_m^TB B^T\bV_m Y_m    + \bV_m^TM^{-1}\Pi C^TC \Pi^TM^{-1}\bV_m = 0.$$
When we apply  the extended block Arnoldi Algorithm \ref{thetaalgo} to the pair $(M_{\Theta}^{-1}A_{\Theta}^T,M_{\Theta}^{-1}C_{\Theta}^T)$, we can notice that $M_{\Theta}^{-1} C_{\Theta}^T=\cV_1^{(1)} \Lambda^{(1,1)}$ resulting from the use of {\tt qr} function in Step 2 and then
\begin{align*}
	\Theta_r M_{\Theta}^{-1} C_{\Theta}^T&= \Theta_r \cV_1^{(1)} \Lambda^{(1,1)}, \\
	\Theta_r M_{\Theta}^{-1} \Theta_r^T C^T&= V_1^{(1)} \Lambda^{(1,1)}.
\end{align*}
We  already proved that $\Theta_r M_{\Theta}^{-1} \Theta_r^T= M^{-1} \Pi$, which gives 
\begin{align}
	M^{-1} \Pi\, C^T&= V_1^{(1)} \Lambda^{(1,1)}, \label{piC}\\
	\bV_m^T M^{-1} \Pi\, C^T&=\bV_m^T V_1^{(1)} \Lambda^{(1,1)}= \begin{bmatrix} I_{n_c}  \\ 0_{n_c} \\ \vdots \\ 0_{n_c}  \end{bmatrix} \Lambda^{(1,1)}. 
\end{align}
Finally, we  obtain  the following low-dimensional Riccati equation
\begin{align}\label{lowdimricca}
	\bT_{m} Y_m  +   Y_m \bT_{m}^T -  Y_m \bV_m^TB B^T\bV_m Y_m + \bV_m^T V_1^{(1)} \Lambda^{(1,1)}C^T (\bV_m^T V_1^{(1)} \Lambda^{(1,1)}C)^T=0,
\end{align}
which is solved by a direct method such as {\tt care} in MATLAB. \\
Let  $R(X_m)$ be the residual corresponding to the approximation $X_m$ given by
\begin{subequations}
	\label{residual}
	\begin{eqnarray}
		R(X_m) &=& M^{-1}\Pi A^T X_m +  X_m A \Pi^T M^{-1}\\ 
		&-&X_mB B^TX_m   + M^{-1}\Pi C^TC \Pi^TM^{-1}.
	\end{eqnarray}
\end{subequations}

\noindent In order to stop the iterations, we need to compute the residual $R(X_m)$ given by (\ref{residual}) without involving $X_m$, since it becomes expensive as $m$ increases. The next result shows how to compute the residual norm of $R(X_m)$ without involving the approximate solution,  which is given only in a factored form at the end of the process.

\begin{theorem1}
	Let $\bV_m \in \R^{2mn_c \times 2mn_c}$ be an orthonormal matrix generated by Algorithm\ref{newextd}. Let $X_m=\bV_m Y_m \bV_m^T$ be the approximate solution of the GARe (\ref{riccatinew}), then the residual norm is given by
	\begin{equation}
		\|R(X_m)\|=\|T_{m+1,m} E_m^T Y_m\|,
	\end{equation}
	where $E_m = [0_{2n_c \times 2(m-1)n_c},I_{2n_c}]^T$ and $\|\|$ is the abbreviation of $\|\|_2.$ 
\end{theorem1}

\begin{proof}
	According to  (\ref{equaalgTm}) and (\ref{residual}), we have 
	\begin{eqnarray*}
		R(X_m)&=&  M^{-1}\Pi A^T X_m +  X_m A \Pi^T M^{-1}-X_mB B^TX_m +M^{-1}\Pi C^TC \Pi^TM^{-1} \\
		&=& M^{-1}\Pi A^T \bV_m Y_m \bV_m^T +  \bV_m Y_m \bV_m^T A \Pi^T M^{-1} - \bV_m Y_m \bV_m^TB B^T\bV_m Y_m \bV_m^T   + M^{-1}\Pi C^TC \Pi^TM^{-1} \\
		&=& (\bV_m \bT_m +V_{m+1} T_{m+1,m} E_m^T)Y_m \bV_m^T + \bV_m Y_m (\bT_m^T \bV_m^T + E_m T^T_{m+1,m} V_{m+1}^T) \\
		&-& \bV_m Y_m \bV_m^TB B^T\bV_m Y_m \bV_m^T 
		+M^{-1}\Pi C^TC \Pi^TM^{-1}.\\
	\end{eqnarray*}
	Using the fact that $ M^{-1} \Pi\, C= V_1^{(1)} \Lambda^{(1,1)}$ as it is described in (\ref{piC}), we get
	\begin{eqnarray*}
		R(X_m)&=& [\bV_m, V_{m+1}] \begin{bmatrix}
			\bT_m Y_m  +  Y_m \bT_m^T + \, E_1  \Lambda^{(1,1)}  (E_1  \Lambda^{(1,1)})^T	& (T_{m+1,m} E_m^T Y_m)^T \\ 
			T_{m+1,m} E_m^T Y_m	& 0 
		\end{bmatrix}  \begin{bmatrix}
			\bV_m\\ 
			V_{m+1}
		\end{bmatrix}. 
	\end{eqnarray*}
	Since  $Y_m$ is the symmetric solution of the low-dimensional Riccati equation, then
	\begin{eqnarray*}
		R(X_m)&=& \bV_{m+1}  \begin{bmatrix}
			0	& (T_{m+1,m} E_m^T Y_m)^T \\ 
			T_{m+1,m} E_m^T Y_m	& 0 
		\end{bmatrix} \bV_{m+1}^T,
	\end{eqnarray*}
	and finally we get the desired result 
	\begin{align}
		\label{reresi}
		\|R(X_m)\| = \|T_{m+1,m} E_m^T Y_m\|.
	\end{align}
\end{proof}
We can check weather we get the desired convergence by verifying the test $\|R(X_m)\|< \epsilon$. Fortunately, the residual $\|R(X_m)\|$ can be computed in a suitable way as described in the theorem above, without computing the approximate solution $X_m$. We take the advantage of $X_m$ as a symmetric positive semi-definite, so it can be decomposed into a  product of two matrices of low-rank as $X_m=Z Z^T$, where $Z$ is a matrix of rank less than or equal to $2m$. The benefit from this decomposition is that we just need to store $Z$ in order to compute the approximate solution $X_m$. Let $Y_m= U \Sigma V$, the SVD decomposition of $Y_m$ where $\Sigma$ is the  matrix of the singular values of $Y_m$ sorted in decreasing order. Let  {\tt dtol} some tolerance and  define $U_r$, $ V_r$ as the first $r$ columns respectively of  $U$ and $V$ corresponding to the $r$ singular values of magnitude greater than {\tt dtol}. In the numerical experiments, we set dtol=$10^{-12}$. Setting $\displaystyle{\Sigma_r}=[\sigma_1, \cdots, \sigma_r]$, we get  $Y_m \approx U_r \Sigma_r V_r^T$, and it follows that
\begin{align}
	\label{ZZT}
	X_m \approx Z_m Z_m^T,
\end{align}
with $Z_m= \bV_m U_r (\Sigma_r)^{1/2}.$ 

\noindent The iterations were stopped when the relative residual norm was less than {\it tol}$\, =10^{-8}$
\begin{equation}
	\label{ter}
	\frac{\|R(X_m)\|}{\|M^{-1}\Pi C^T C\Pi^TM^{-1}\|} < 10^{-8}.
\end{equation}
We mentioned before that all our results are obtained without any explicit computation of $\Pi$, so to compute $M^{-1}\Pi C^T$ in an appropriate manner we use the formula (\ref{piC}), and then 
\begin{align}
	\label{MPic}
	M^{-1}\Pi C^T C \Pi^T M^{-1} = V_1^{(1)} \Lambda^{(1,1)} (V_1^{(1)} \Lambda^{(1,1)})^T,
\end{align}
where $ V_1^{(1)}$ is the first $n_c$ block-column of $V_1 \in \R^{n_v \times 2n_c}$ and $\Lambda^{(1,1)} \in \R^{n_c \times n_c}$ is the block $(1,1)$ of the upper triangular matrix $\Lambda \in \R^{2n_c \times 2n_c}$ previously described in (\ref{step1algo}).   
The following algorithm summarizes all the results explained above.

\begin{algorithm}[H]
	\caption{ Extended block Arnoldi Riccati algorithm (EBARA)}
	\label{ERBALA}
	\begin{itemize}		
		\item[$\bullet$]  Inputs: ~$M,\, A\in \R^{n_v \times n_v}$, $G \in \R^{n_v \times n_p}$, $B\in \R^{n_v \times n_b}$, $C \in \R^{n_c \times n_v}$, tolerance $\epsilon$, dtol, number of iteration $m_{max}$.
		\item[$\bullet$]  Outputs: the approximate solution $X_m \approx Z_m Z_m^T$.
		\item[$\bullet$] For $m=1, \cdots, m_{max}$ 
		\item[$\bullet$] Use Algorithm \ref{newextd} to compute $\bV_m$ an orthonormal matrix and compute $\bT_m$ the block Hessenberg matrix.
		\item[$\bullet$] Solve the low-dimensional Riccati equation (\ref{lowdimricca}) using the MATLAB function {\tt care}.
		\item[$\bullet$] Compute the relative residual norm (\ref{ter}) using  (\ref{reresi}) and (\ref{MPic}), and if it is less than $\epsilon$, then
		\begin{enumerate}
			\item Compute the SVD of $Y_m= U \Sigma V$ where $\Sigma= diag [\sigma_1, \cdots, \sigma_{2m}].$
			\item Determine $r$ such that $\sigma_{r+1} < \text{dtol} \leq \sigma_r$, set $\Sigma_r =diag[\sigma_1, \cdots, \sigma_r]$ and compute $Z_m= \bV_m U_r (\Sigma_r)^{1/2}$.
		\end{enumerate} 
		end if.
		\item[$\bullet$] End For
	\end{itemize}
	\label{ebaraalgo}
\end{algorithm}
\section{Numerical experiments}
In this section, we present some  numerical results to confirm the performance of the proposed  approaches. All the experiments were carried out using MATLAB R2018a on a computer with Intel $^\text{\textregistered}$ core i7 at 2.3GHz and 8Gb of RAM. The MATLAB programs representing the two algorithms (Algorithm \ref{newextd}, Algorithm \ref{ERBALA}) are available in \url{https://lmpa.univ-littoral.fr/index.php?page_id=8}. Our method is applied to a discretized Navier-Stokes equations as described in Section \ref{sec2}.  We first show how our method allows us to build an efficient reduced model by presenting the transfer functions of the original and reduced systems with the associated error. Then we  investigate the numerical solution of the GARe (\ref{riccatinew}) using  Algorithm \ref{ebaraalgo} and as we mentioned earlier  this numerical solution is actually the key to construct the matrix feedback used to stabilize the unstable  system.  All the data was provided from \cite{benner15}. Some information on this data are depicted in Table \ref{tabnvnp}. The state dimension $n_v$ refers to the dimension of the discretized velocity field, and $n_p$ is the dimension of the discretized pressure field, also sparsity of each matrix $A$ and $M$ is given, i.e., the ratio of the number of non-zero elements to the total number of elements in the matrix.
We used different dimensions of $n_v$ and $n_p$ corresponding to three levels. 
\begin{table}[!h]
	\begin{center}
		\caption{The matrix dimensions for different levels}
		\label{tabnvnp}
		\begin{tabular}{c|c|c|c|c}
			\hline
			Level 	& $n_v$ & $n_p$ & full model ($n_v+n_p$) & sparsity of  $A \, \& \, M$ \\
			\hline 
			1	& 4796 & 672 & 5468 & $4.6\cdot 10^{-3} \,  \mid \, 2.3 \cdot 10^{-3} $ \\
			2	& 12292 & 1650 & 13942 &$1.8 \cdot 10^{-3}\,  \mid \, 9.05\cdot10^{-4}$ \\
			3  & 28914 & 3784 & 32698 & 7.79$\cdot10^{-4}\,  \mid \,3.89\cdot10^{-4}$\\
			\hline
		\end{tabular} 
	\end{center}
\end{table}
Note that the norm used here is the $\cH_{\infty}$ norm and it expressed as $\| F- F_m\|_{\infty}= \displaystyle\sup_{\omega \in \R}\|F(j\omega)-F_m(j\omega)\|_2$. To compute this norm we use the following functions from \texttt{lyapack} \cite{laypack} 
\begin{enumerate}
	\item {\tt lp$\_$lgfrq} : Generates a set of logarithmically distributed frequency sampling points $\omega \in [10^{-5},10^5]$.
	\item {\tt lp$\_$gnorm} : Computes a vector which contains the 2-norm $$\| F-F_m\|= \sigma_{max}(F(i\omega)-F_m(i\omega)),$$ for each sampling points $\omega\in [10^{-5},10^5]$ and $i=\sqrt{-1}$.
\end{enumerate}

\begin{figure}[h]
	\centering
	\includegraphics[width=13cm]{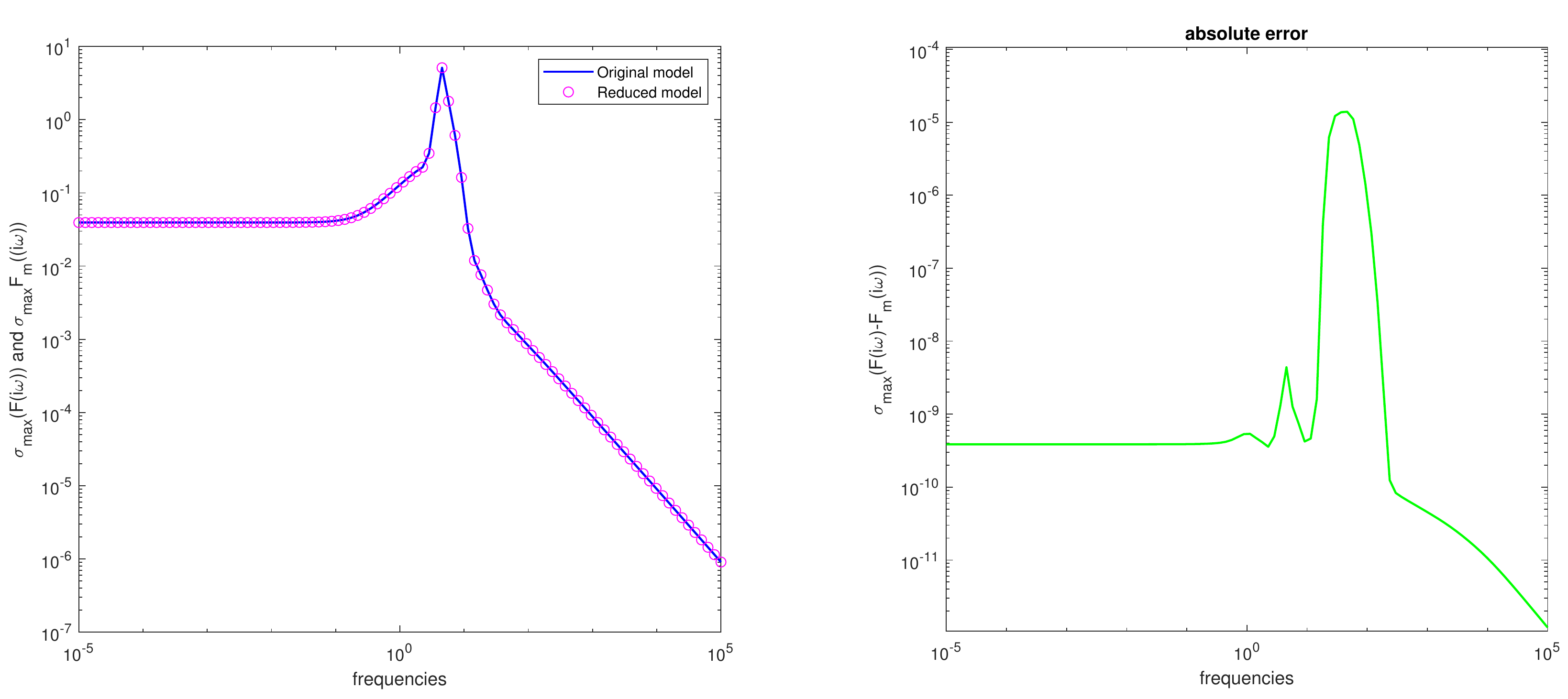}
	\centering
	\caption{Level 1 with Re= 300: Bode plot (left) and the error norms versus frequencies (right).}
	\label{fig1}
\end{figure}

\begin{figure}
	\centering
	\includegraphics[width=13cm,height=6.5cm]{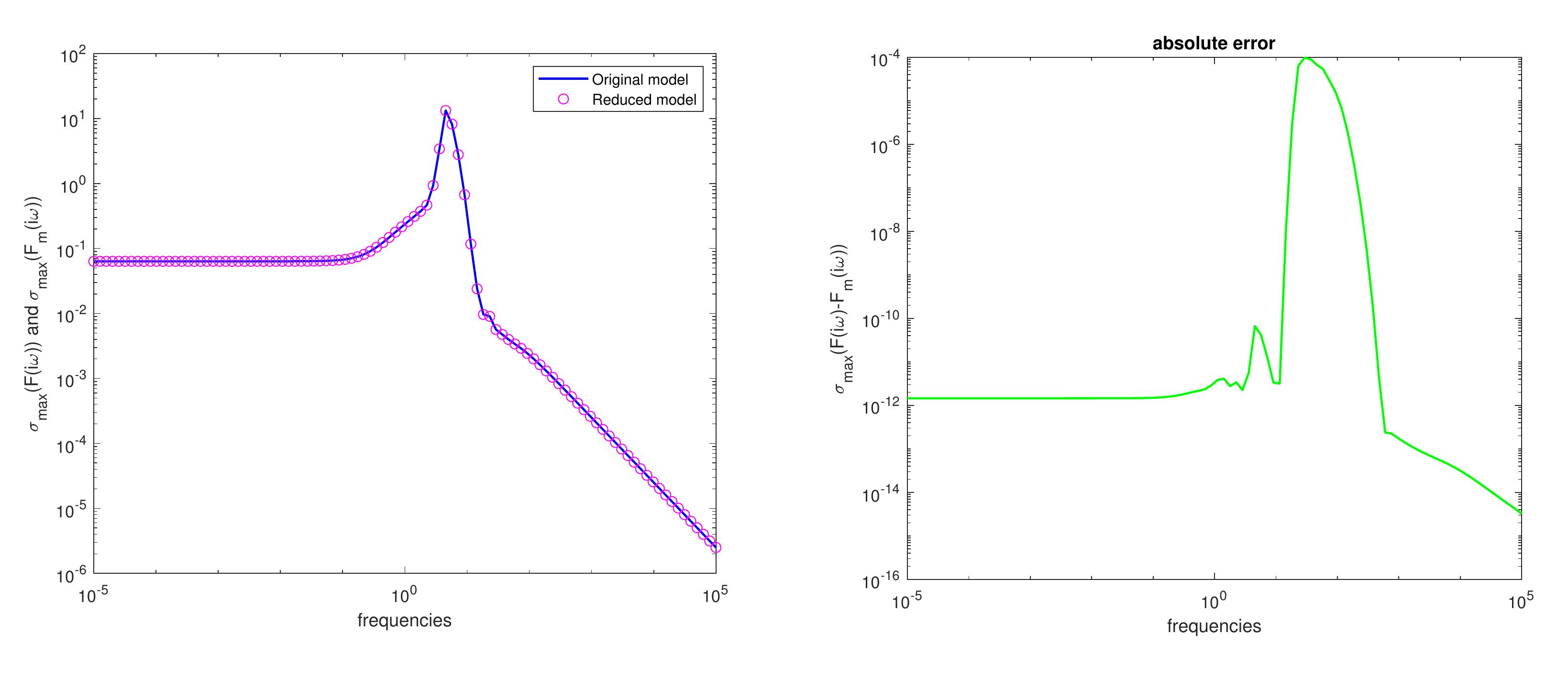}
	\caption{Level 2 with Re= 400: Bode plot (left) and the error norms versus frequencies (right). }
	\label{fig2}
\end{figure}

\begin{figure}
	\centering
	\includegraphics[width=13cm]{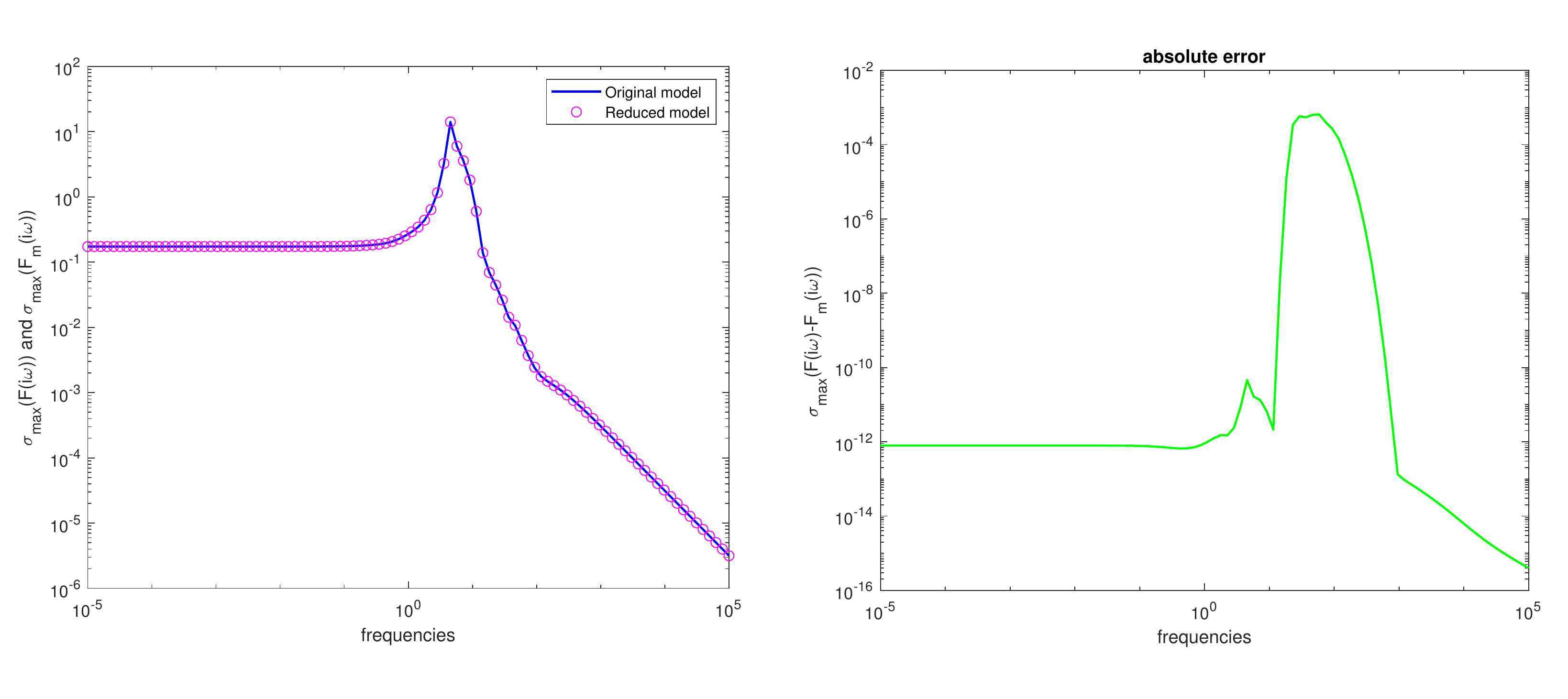}
	\caption{Level 3 with Re= 500: Bode plot (left) and the error norms versus frequencies (right).}
	\label{fig3}
\end{figure}

{\bf Example 1} For this example, we show the frequency response of the original and reduced systems. We considered  the three models from Table \ref{tabnvnp} where we associate level 1 with Reynolds number $\text{Re}=300$, level 2 with $\text{Re}=400$ and level 3 with  $\text{Re}=500$. For the three models we used  $m=120$ and  then the dimension of the reduced system is $2\times m\times n_b=480$. For a Reynolds number $\text{Re}<100$,  Navier-Stokes flow starts to behave like a Stokes flow, and this comes from the fact that the convection term $(z \cdot \nabla )z$ in (\ref{NSequa}) doesn't have an important impact. Figures \ref{fig1}, \ref{fig2} and \ref{fig3}  illustrate  the obtained results comparing the original transfer function and its approximation for the three levels.  We also plotted the error-norm  between the two transfer functions. The computed error norm  $\| F- F_m\|_{\infty}$ was  $1.37 \cdot 10^{-5}$ for level 1, $\| F- F_m\|_{\infty}=9.82 \cdot 10^{-5}$ for level 2 and $\| F- F_m\|_{\infty}= 6.5 \cdot 10^{-4}$ for level 3. 


\begin{figure}[h!]
		\centering
		\includegraphics[width=7cm]{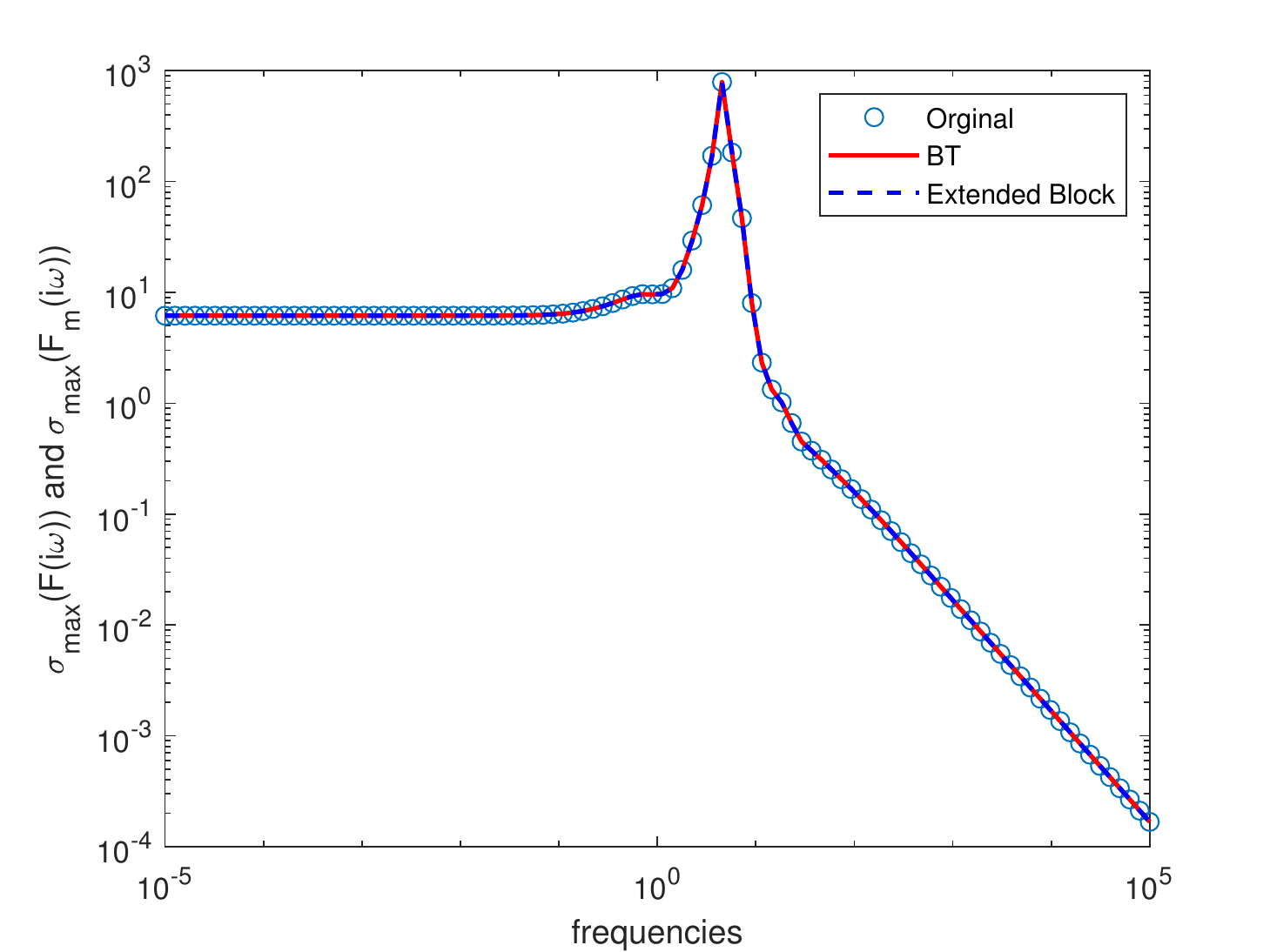}
		\centering
		\caption{Level 2 with Re= 300: Bode plot}
		\label{compfig1}
	\end{figure}

\begin{figure}[h!]
		\includegraphics[width=6.5cm]{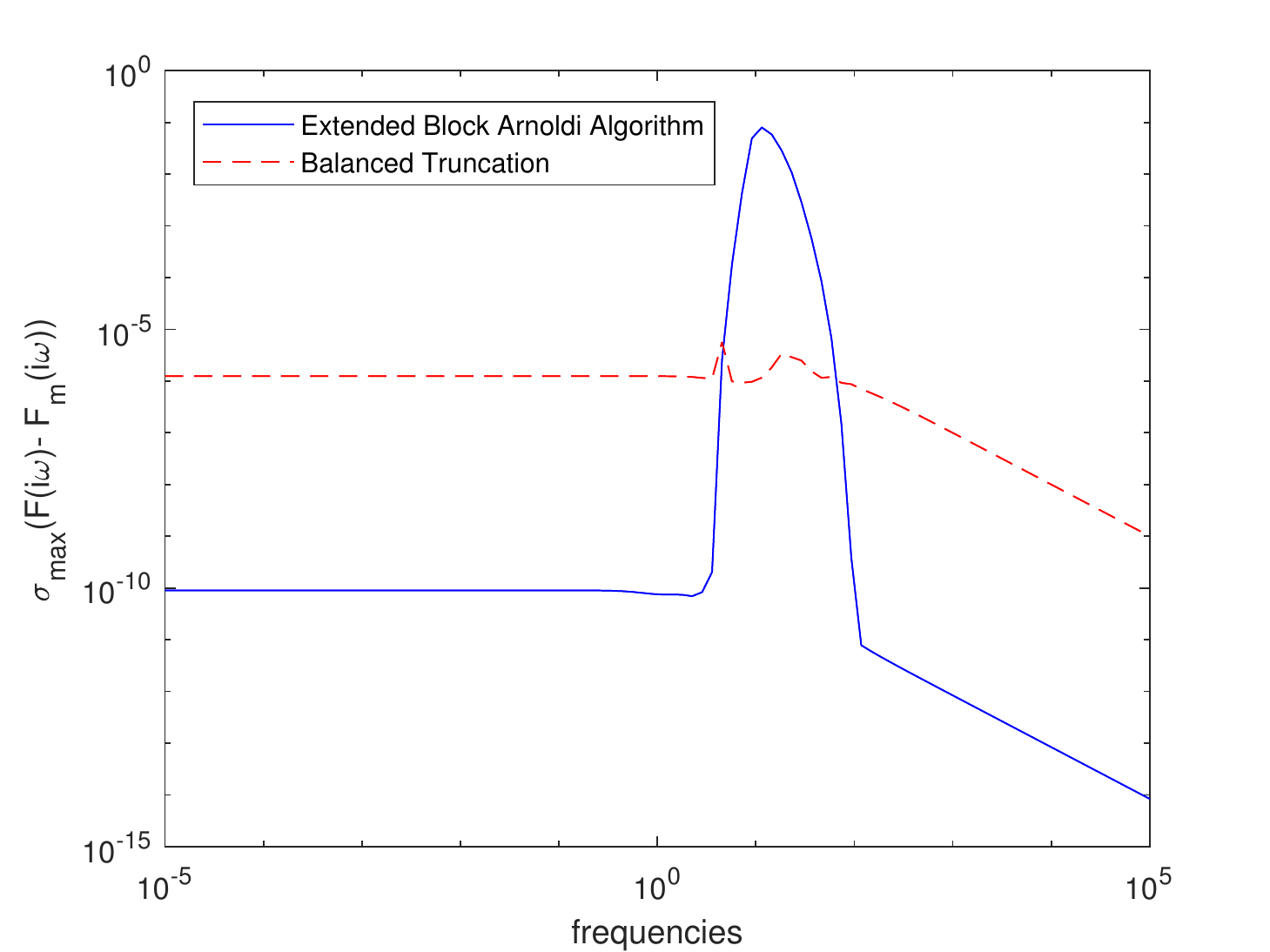}
		\hfill
		\includegraphics[width=6.5cm]{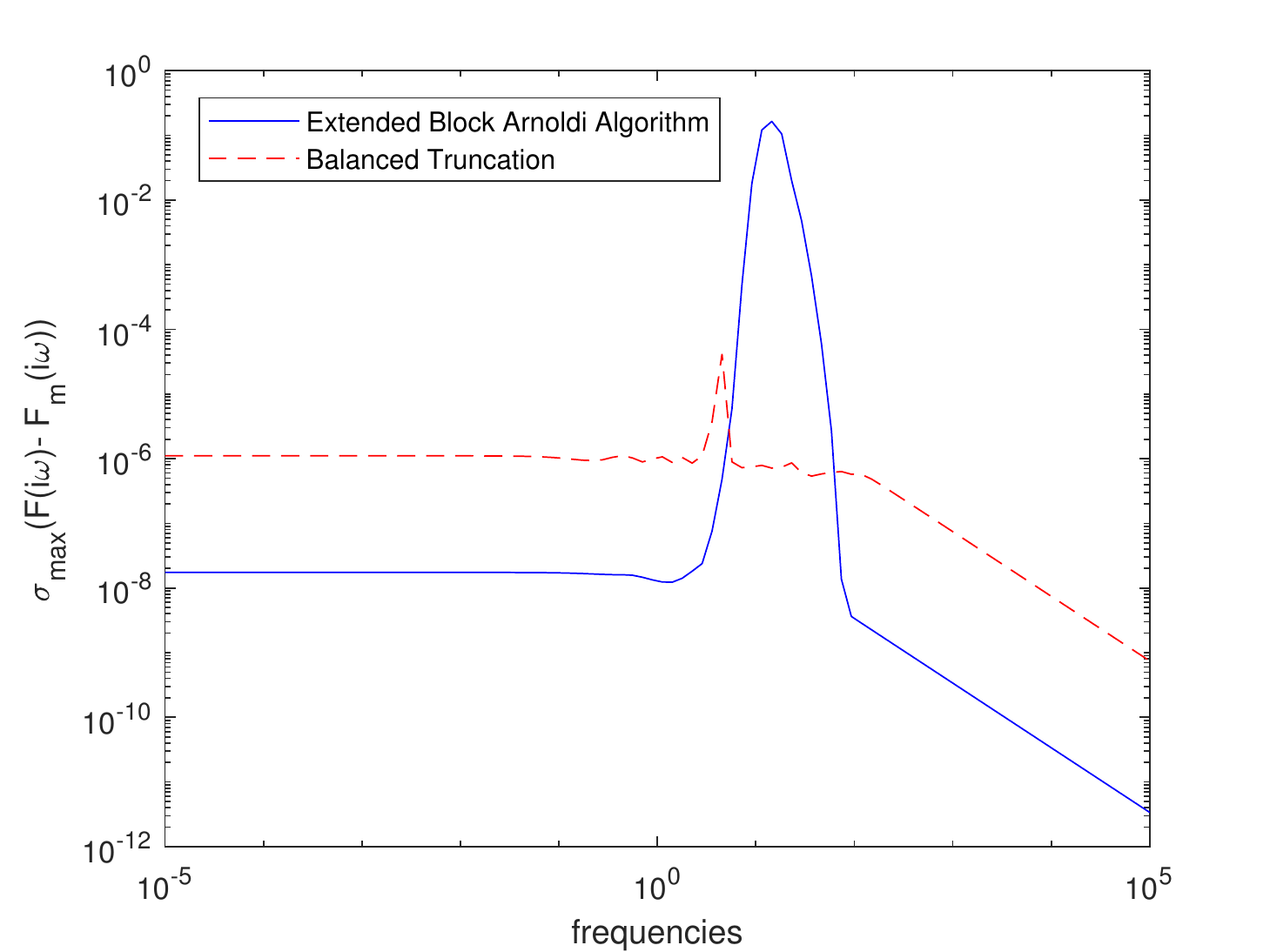}
		\caption{The error-norms versus frequencies using level 1 with \text{Re}=300 (left)  and with Re=400 (right).}
		\label{comp1}
	\end{figure}

\begin{figure}[h!]
		\includegraphics[width=6.5cm]{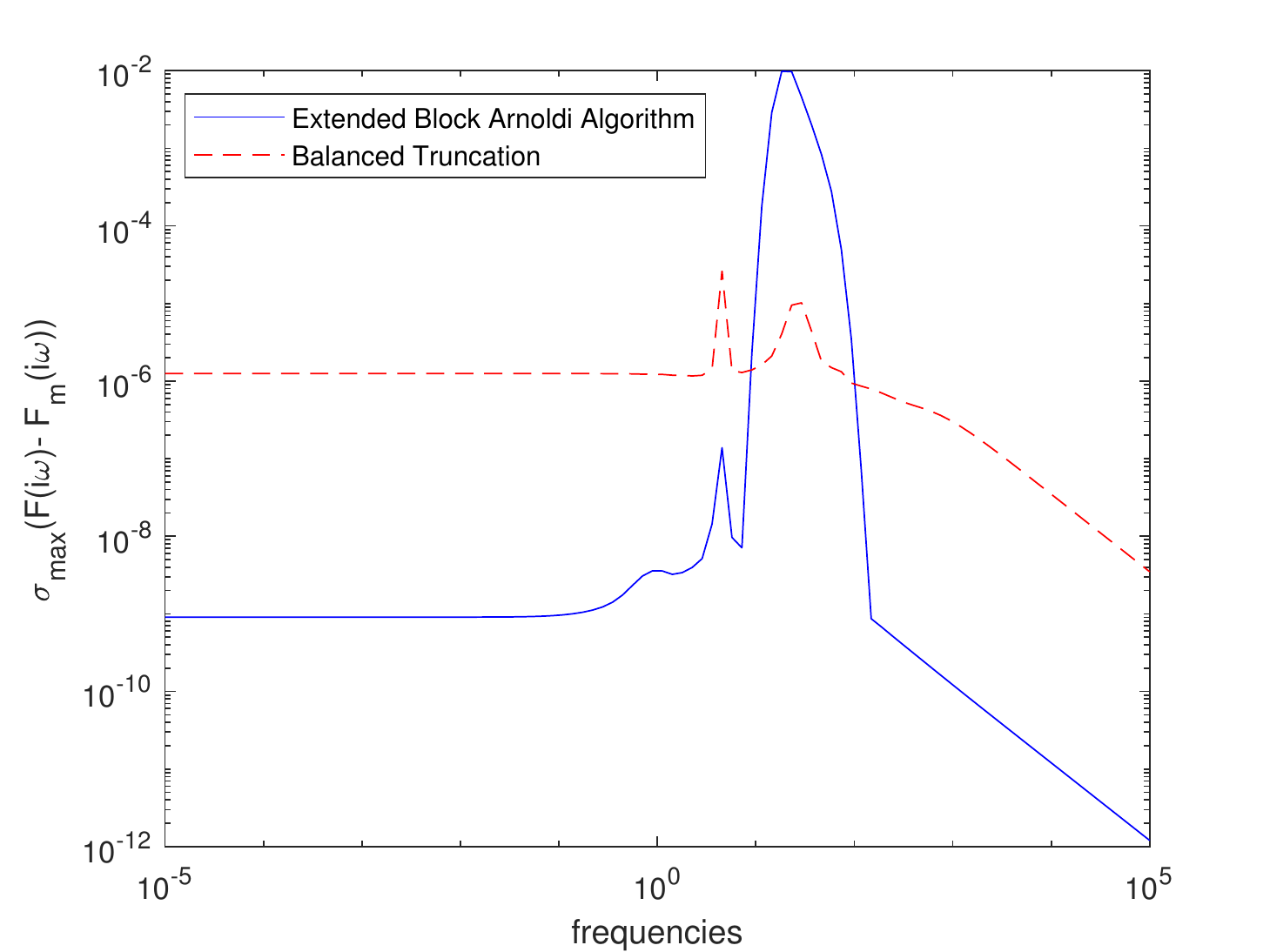}
		\hfill
		\includegraphics[width=6.5cm]{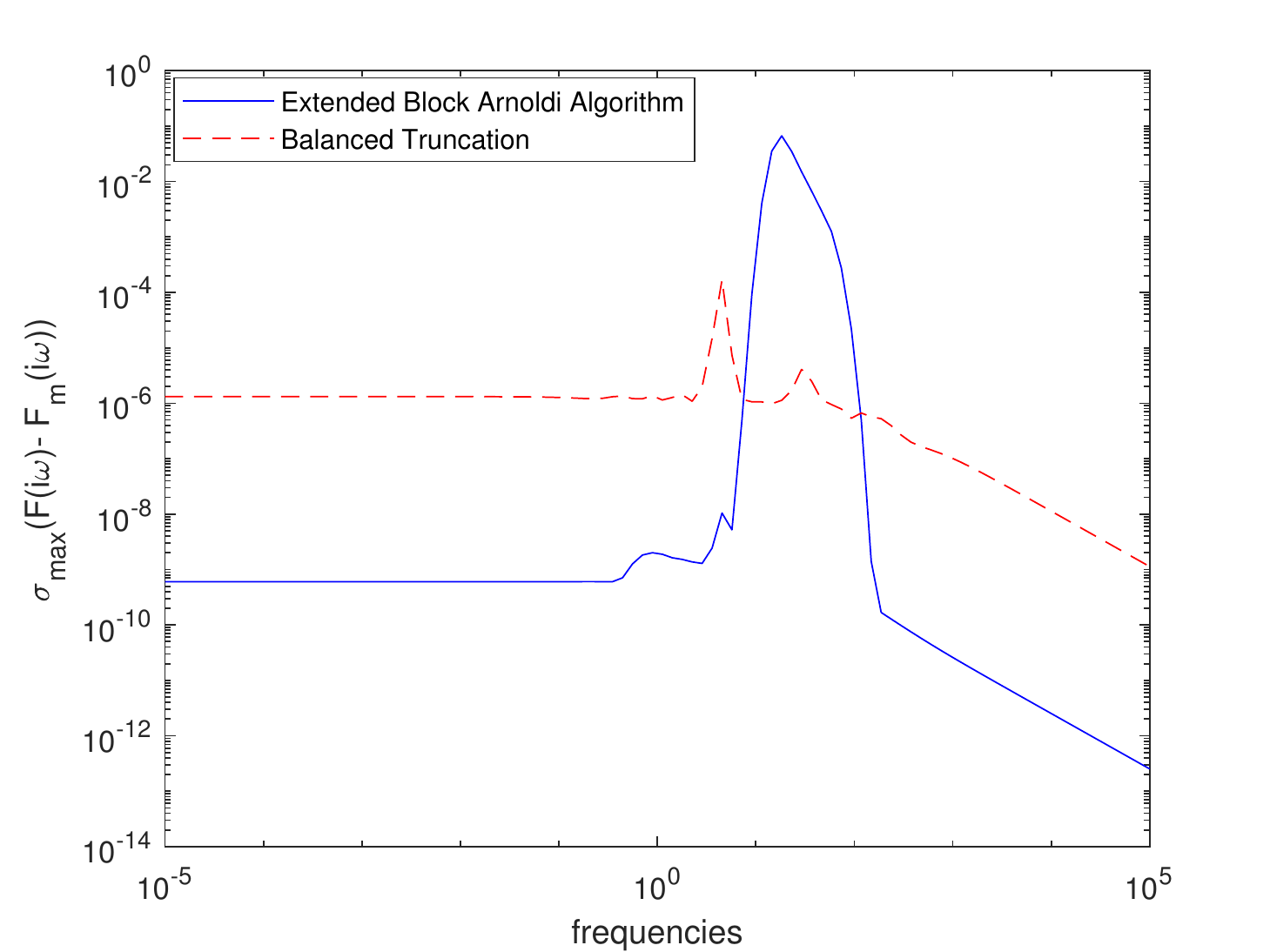}
		\caption{The error-norms versus frequencies using level 2 with \text{Re}=300 (left)  and with Re=400 (right).}
		\label{comp2}
	\end{figure}

{\bf Example 2} To investigate the efficiency of Algorithm \ref{newextd}, we compare our method to a common and deployed model order reduction method knows as the Balanced Truncation (BT). The main challenge in the BT is to solve larges-scale Lyapunov equations in order to obtain the system Gramians that will be used to generate a reduced model. The BT algorithm is available at the M-M.E.S.S. toolbox, see \cite{SaaKmess}. The authors used a different data from those presented in Table \ref{tabnvnp}. We chose from their data two level of discretization and we summarize in Table \ref{tabnvnp1} some information. 
	\begin{table}[h!]
		\begin{center}
			\caption{The matrix dimensions for different levels}
			\label{tabnvnp1}
			\vskip0.2cm
			\begin{tabular}{c|c|c|c}
				\hline
				Level 	& $n_v$ & $n_p$ & full model ($n_v+n_p$) \\
				\hline 
				1	& 3142 & 453 & 3595\\  
				2	& 8268 & 1123 & 9391 \\
				\hline
			\end{tabular} 
		\end{center}
	\end{table}
	For the level 1 we used $m=70$ and $m=75$ for level 2. The tolerance truncation is set to $10^{-5}$. In Figure \ref{compfig1}, we plotted the norms $\Vert F(j\omega)\Vert_2$ and its approximation $\Vert F_m(j\omega)\Vert_2$ for different values of the frequency $\omega \in [10^{-5}, 10^5]$ of the two methods (our method and BT). As can be seen, we have obtained a perfect match between the original transfer function and its approximation for both methods. We show the obtained error-norms $\|F(j\omega)-F_m(j\omega)\|_2 = \sigma_{max}(F(j\omega)-F_m(j\omega))$  for different values of the frequency $\omega$ with the Reynold number \text{Re}=300 in Figure \ref{comp1} and with \text{Re}=400 in Figure \ref{comp2}. Here, you can notice that the error of our Algorithm increases rapidly when $\omega \in (1,10^3)$, we tried to alleviate this problem by increasing the number of iterations "m", but this choice increased the computing time and also increased the error when $\omega \in (10^{-5},1)$ from $10^{-10}$ to $10^{-4}$, and this also applies when $\omega \in (10^{3},10^5).$ This is why we stick with the first choice and do not increase the number of iterations 'm'. We present in Table \ref{cputime} the execution time of our algorithm and that based on BT described in \cite{SaaKmess}.
	\begin{table}[h!]
		\begin{center}
			\caption{The CPU-time (in seconds) required for both methods}
			\label{cputime}
			\vskip0.2cm
			\begin{tabular}{c|c|c}
				\hline
				Reynolds number	& Re=300 & Re=400 \\
				& Algorithm \ref{newextd} \; \; BT &  Algorithm \ref{newextd}  \; \;  BT \\
				\hline 
				Level 1	& 4.55 \qquad \;\;\; 13.92 & 4.17 \qquad \;\;\; 15.47\\  
				Level 2 & 19.71 \qquad \;\;\; 51.09& 17.78 \qquad \;\;\; 56.34 \\
				\hline
			\end{tabular} 
		\end{center}
	\end{table}

{\bf Example 3} In this example, we investigate the extended block Arnoldi Riccati algorithm (EBARA, Algorithm \ref{ERBALA}) for solving  generalized algebraic Riccati equations (GARe) (\ref{riccatinew}) which is needed to compute the matrix feedback of our initial problem. We use matrices corresponding to the level 1 of discretization in Table \ref{tabnvnp} with different Reynolds numbers  $\text{Re}=300,400$ and $500$. We have established a comparison between our Algorithm \ref{ERBALA} and Algorithm 2 (a generalized low-rank Cholesky factor Newton method) described in \cite{benner15}. We have summarized in Table \ref{riccacomp} the number of iterations, ADI and Newton iterations as well as the cpu-time needed to reach the convergence of both methods 

\begin{table}[h]
	\renewcommand\arraystretch{1.3}
	\caption{The obtained results of both methods}
	\label{riccacomp}
	\vskip0.2cm
	\begin{tabular}{l | *{5}{>{\centering}p{1.5cm}|}c}
		\hline
		Methods
		& \multicolumn{3}{c|}{Algorithm \ref{ERBALA}} 
		& \multicolumn{3}{c}{Algorithm 2 in \cite{benner15}} \\
		\hline
		
		&$ \#$ of iter. & cpu-time(sec)   &   Rel. res. &   Newton \& ADI iter. &CPU-time(sec)   &   Rel. res.  \\
		\hline
		\text{Re}=300 & 77 &  136.77& $8.43e^{-08}$ & 8 \& 245 & 323.48&  $9.10e^{-09}$\\ 
		\hline
		\text{Re}=400  &94 & 378.82 & $8.36e^{-08}$ & 20 \& 301 & 1124.43& $2.27e^{-07}$ \\
		\hline 
		\text{Re}=500  & 109& 524.27 & $8.14e^{-08}$ & - & $>$1800 & - \\
		\hline
	\end{tabular}
\end{table}
\subsection*{Stabilizing the unstable system}
\begin{figure}[h!]
	\includegraphics[width=6.5cm,height=5cm]{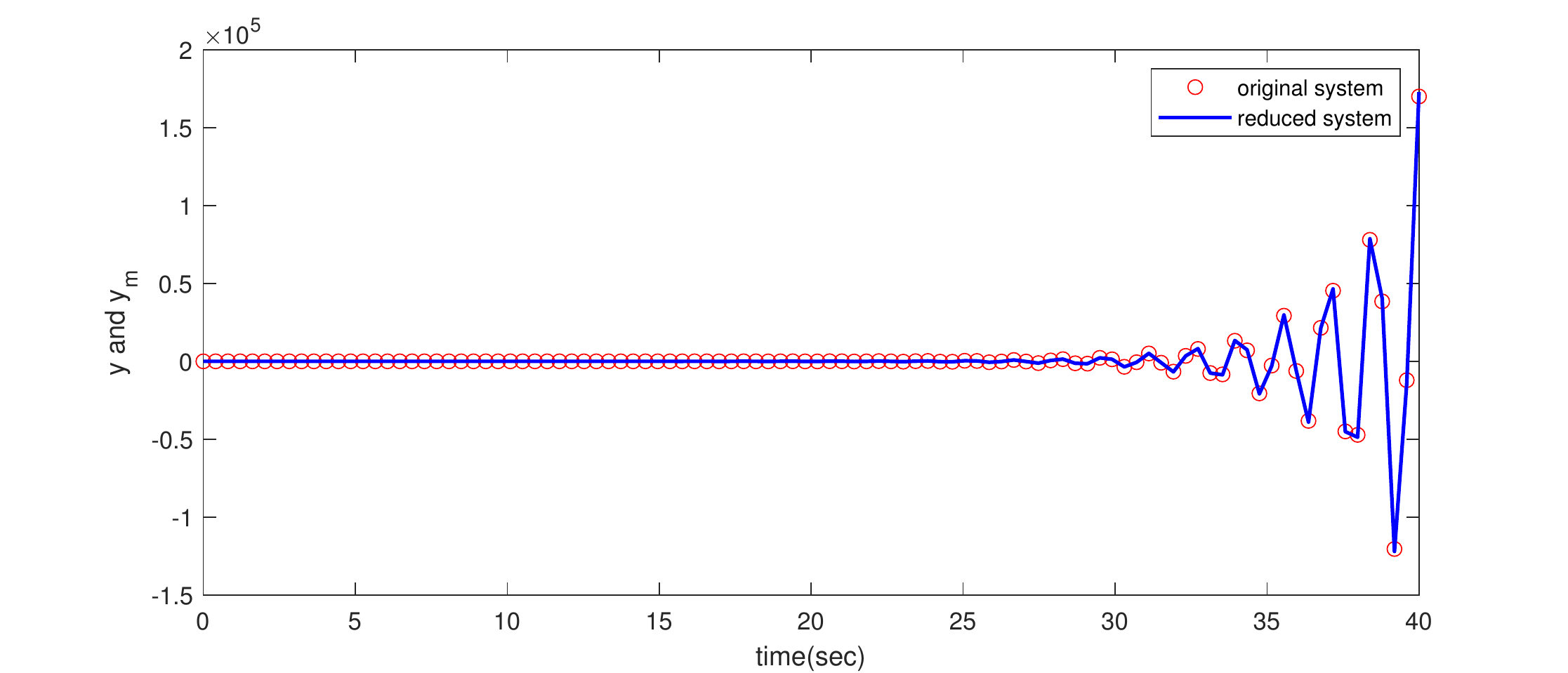}
	\hfill
	\includegraphics[width=6.5cm,height=5cm]{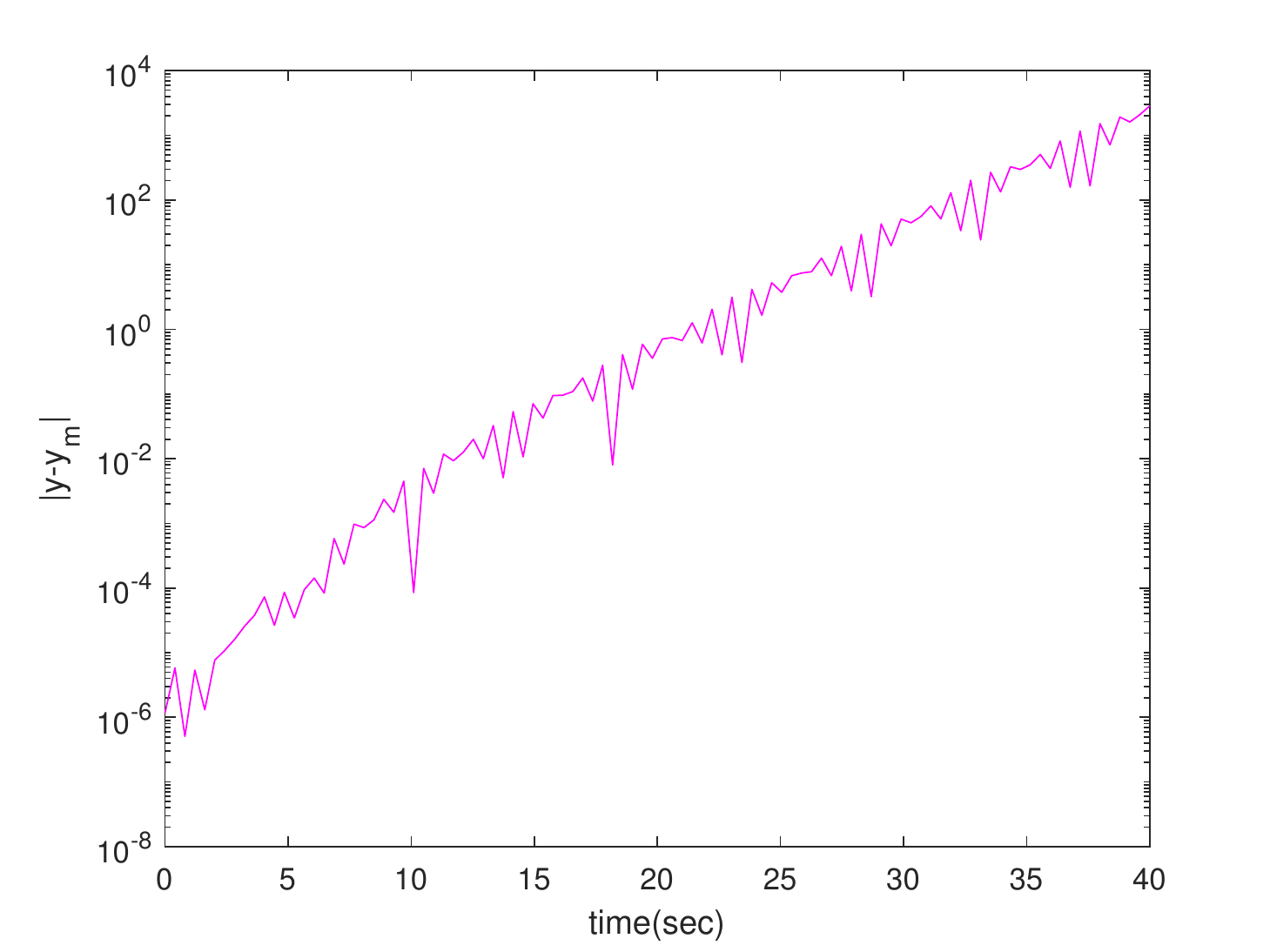}
	\caption{Left: time domain response simulation of the unstable system with Re=400. Right: the error norm $\|y-y_m\|.$}
	\label{bfstab400}
\end{figure}

\begin{figure}[H]
	\includegraphics[width=6.6cm]{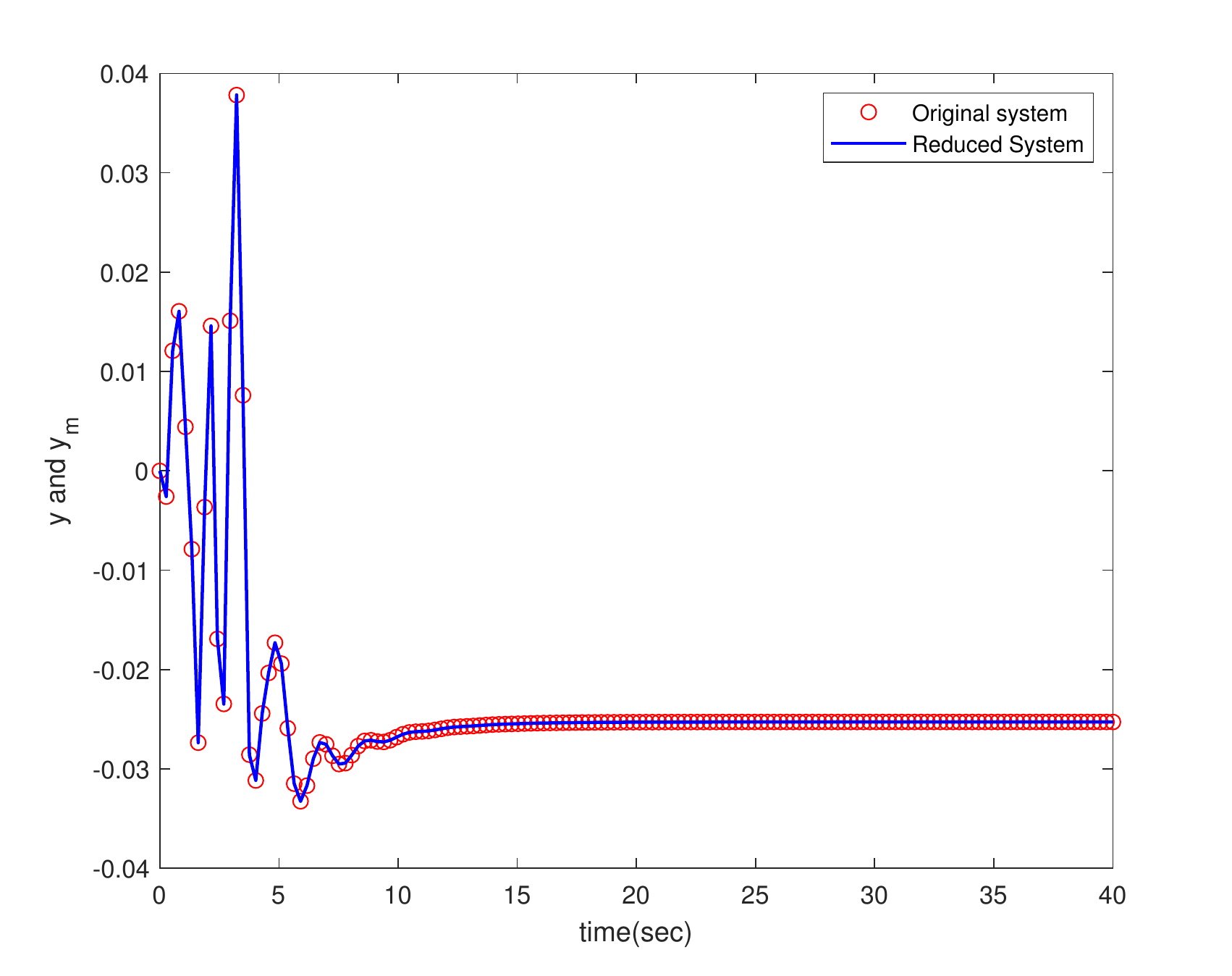}
	\hfill
	\includegraphics[width=6.9cm]{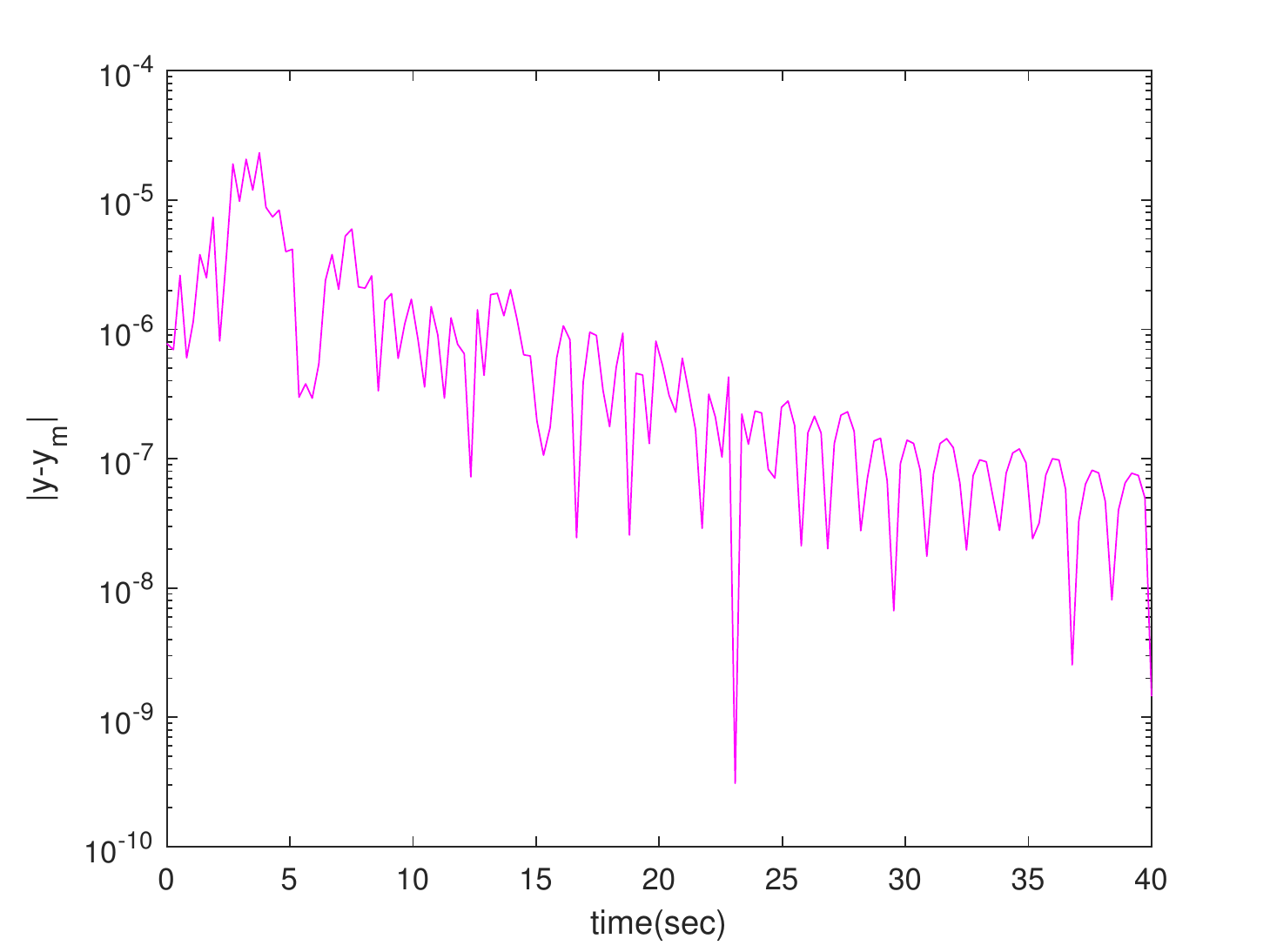}
	\caption{Left: time domain response for the stabilized system (input 1 to output 1) original and reduced systems. Right: the error norm $\|y-y_m\|.$}
	\label{afstab400_1}
\end{figure}

\begin{figure}[h!]
	\includegraphics[width=6.6cm]{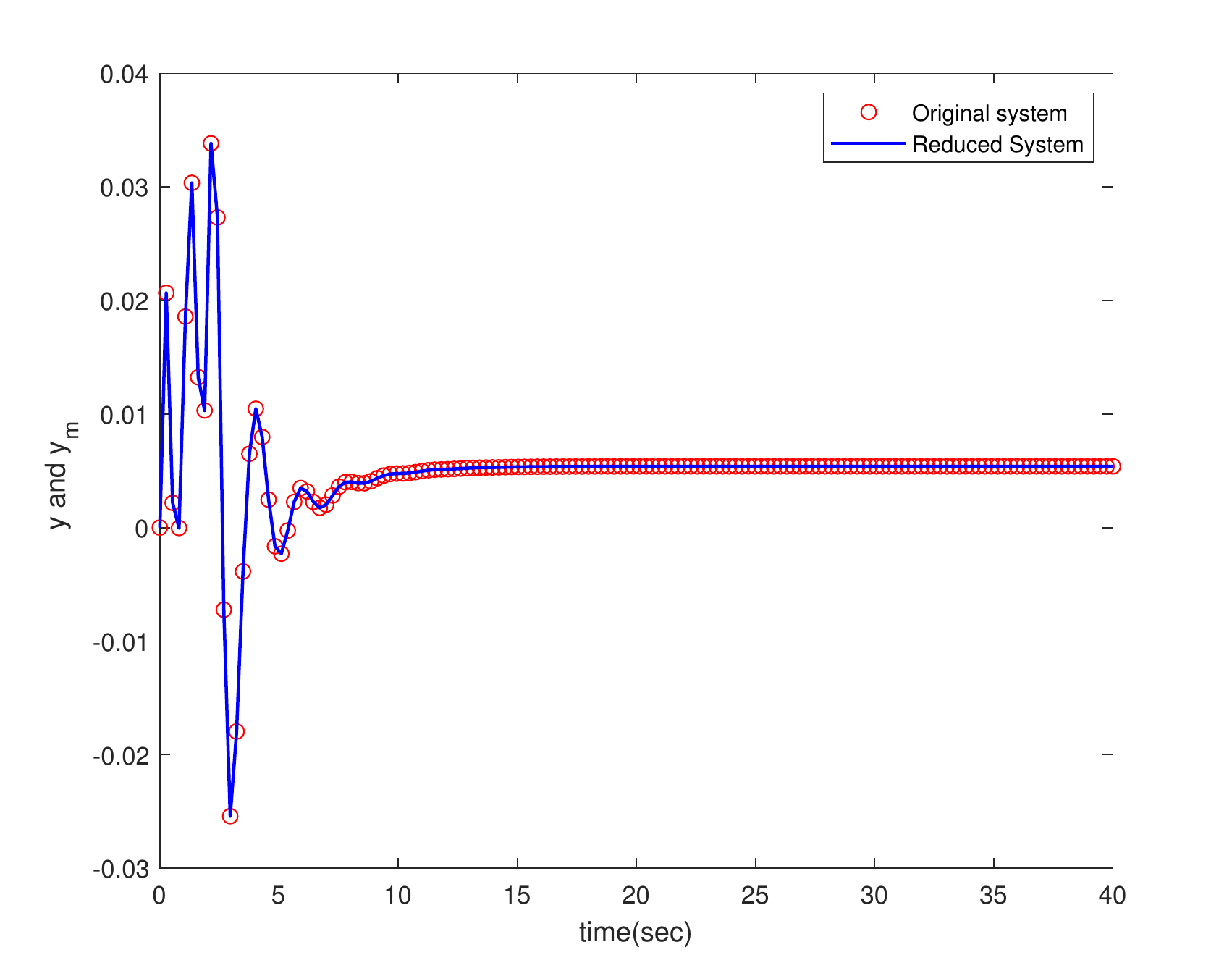}
	\hfill
	\includegraphics[width=6.8cm]{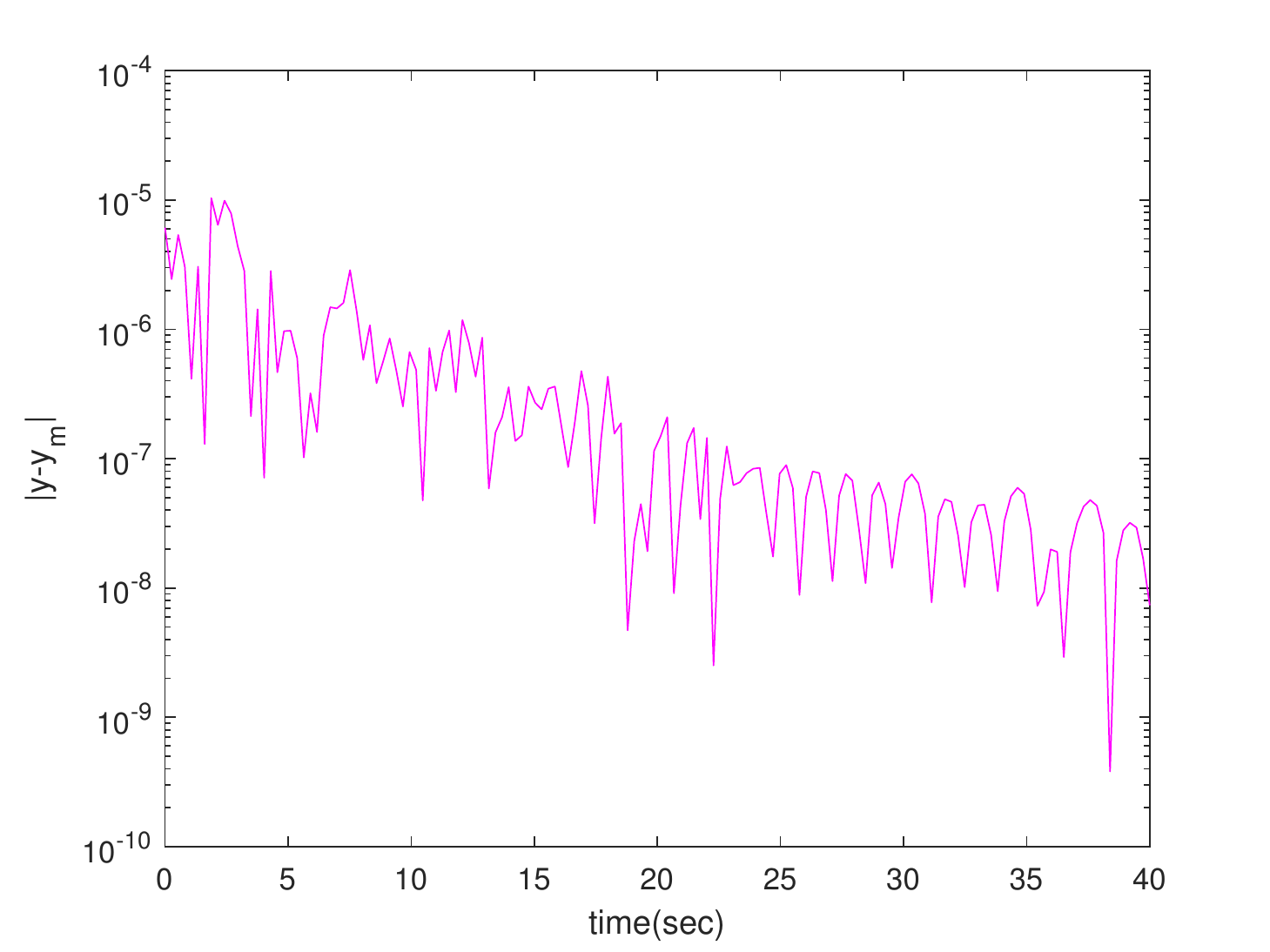}
	\caption{Left: time domain response for the stabilized system (input 2 to output 2) original and reduced systems. Right: the error norm $\|y-y_m\|.$}
	\label{afstab400_2}
\end{figure}

\begin{figure}[h!]
	\includegraphics[width=6.5cm,height=5cm]{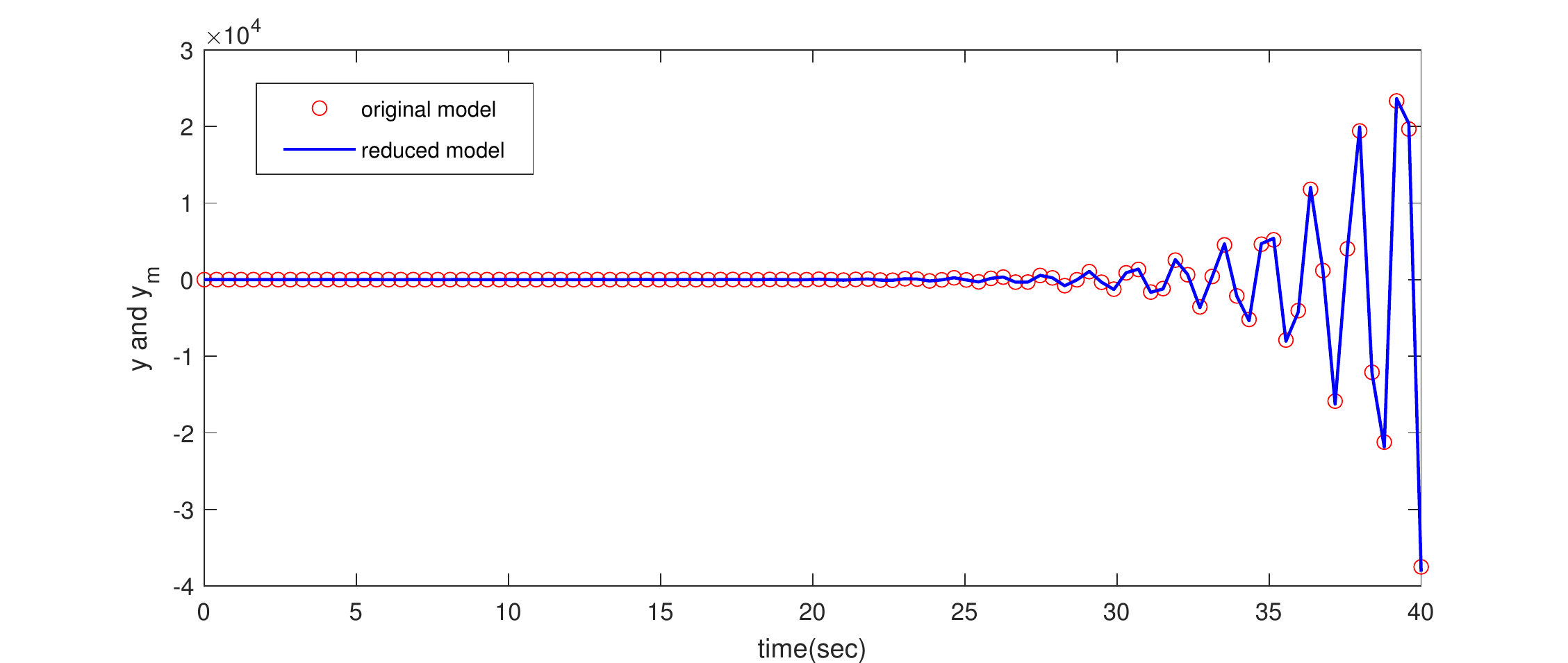}
	\hfill
	\includegraphics[width=6.5cm,height=5cm]{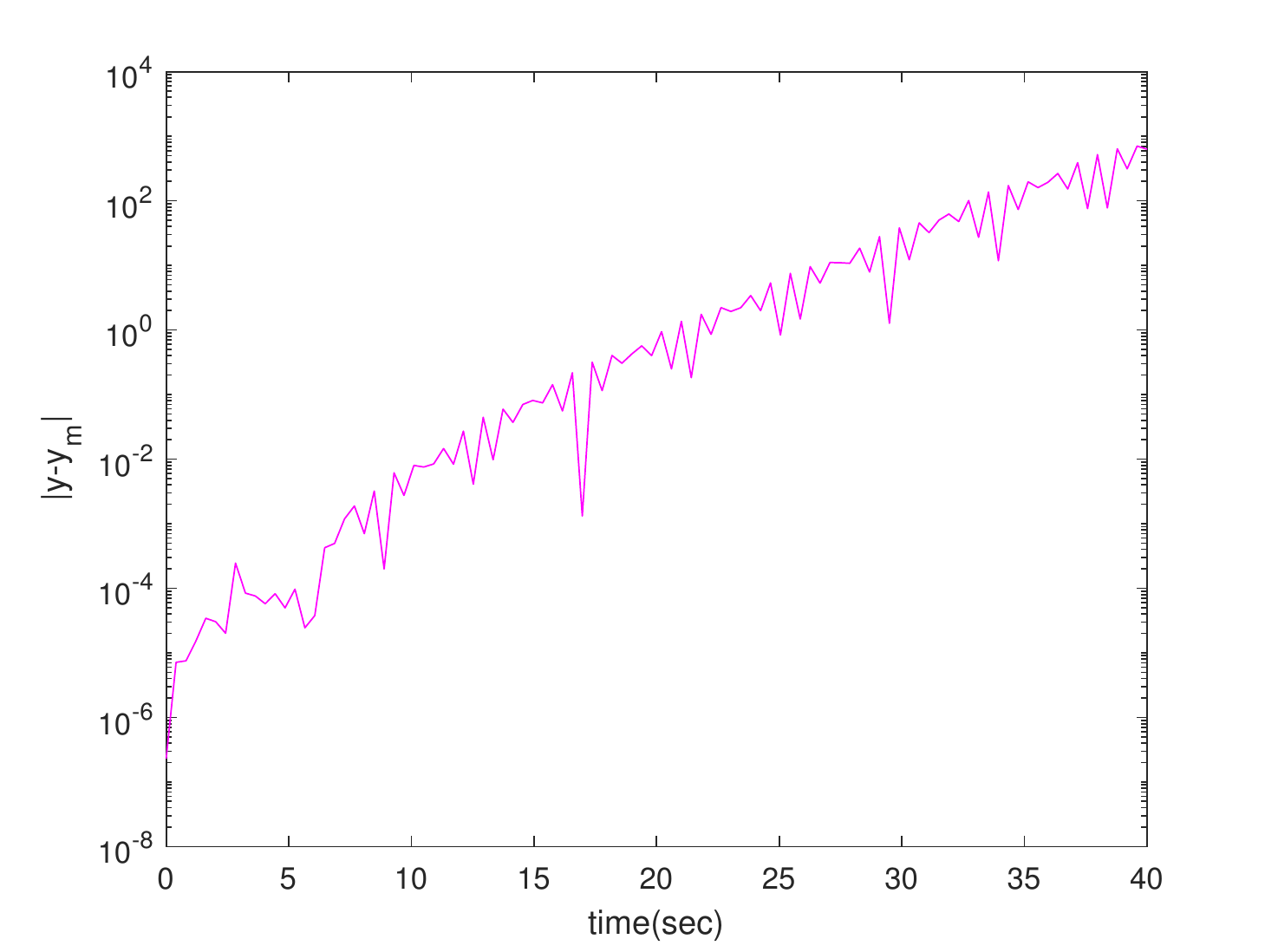}
	\caption{Left: time domain simulation for the unstable system with Re=500. Right: the error norm $\|y-y_m\|.$}
	\label{bfstab500}
\end{figure}

\begin{figure}[H]
	\includegraphics[width=6.5cm,height=4.9cm]{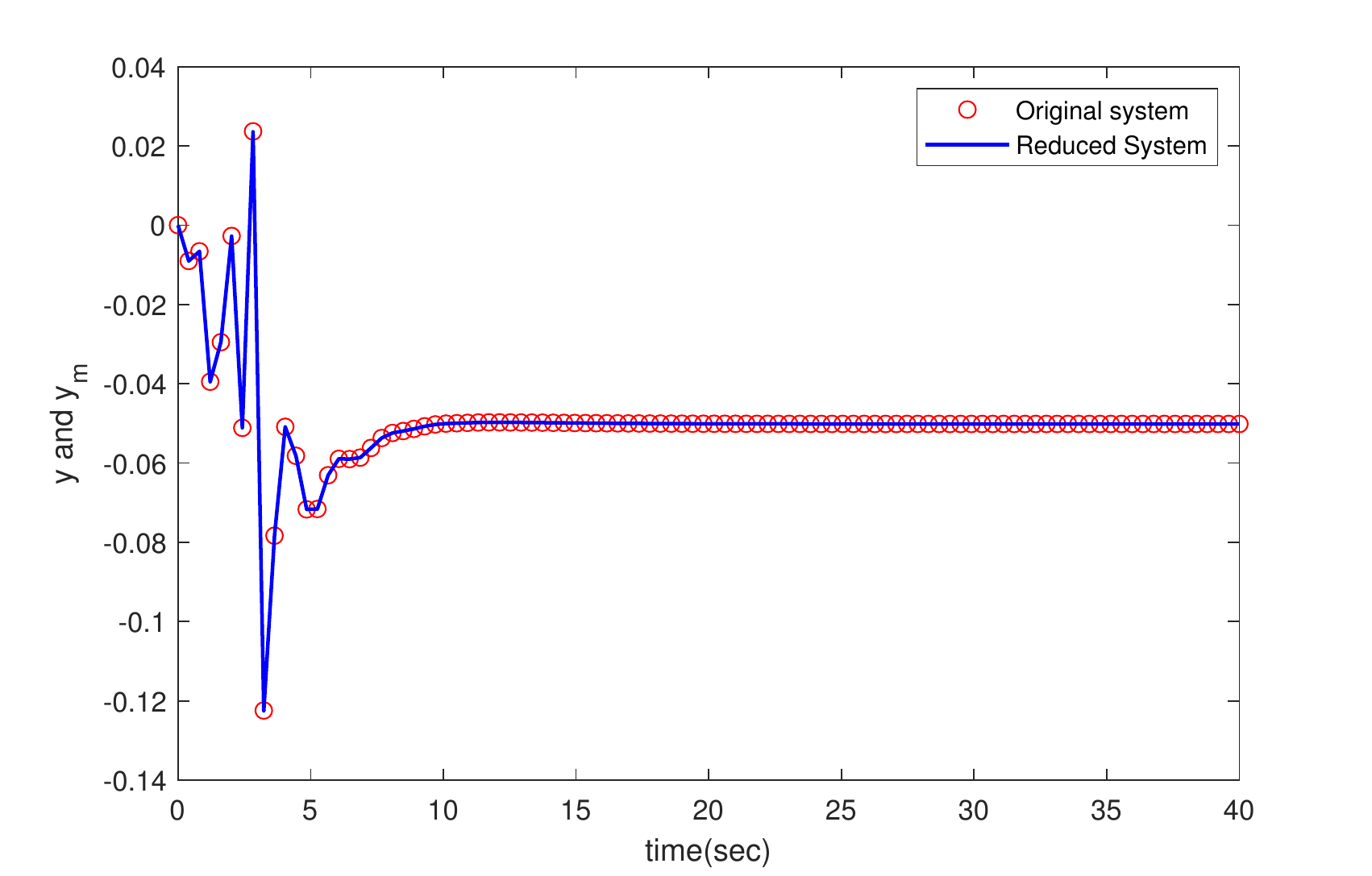}
	\hfill
	\includegraphics[width=6.5cm]{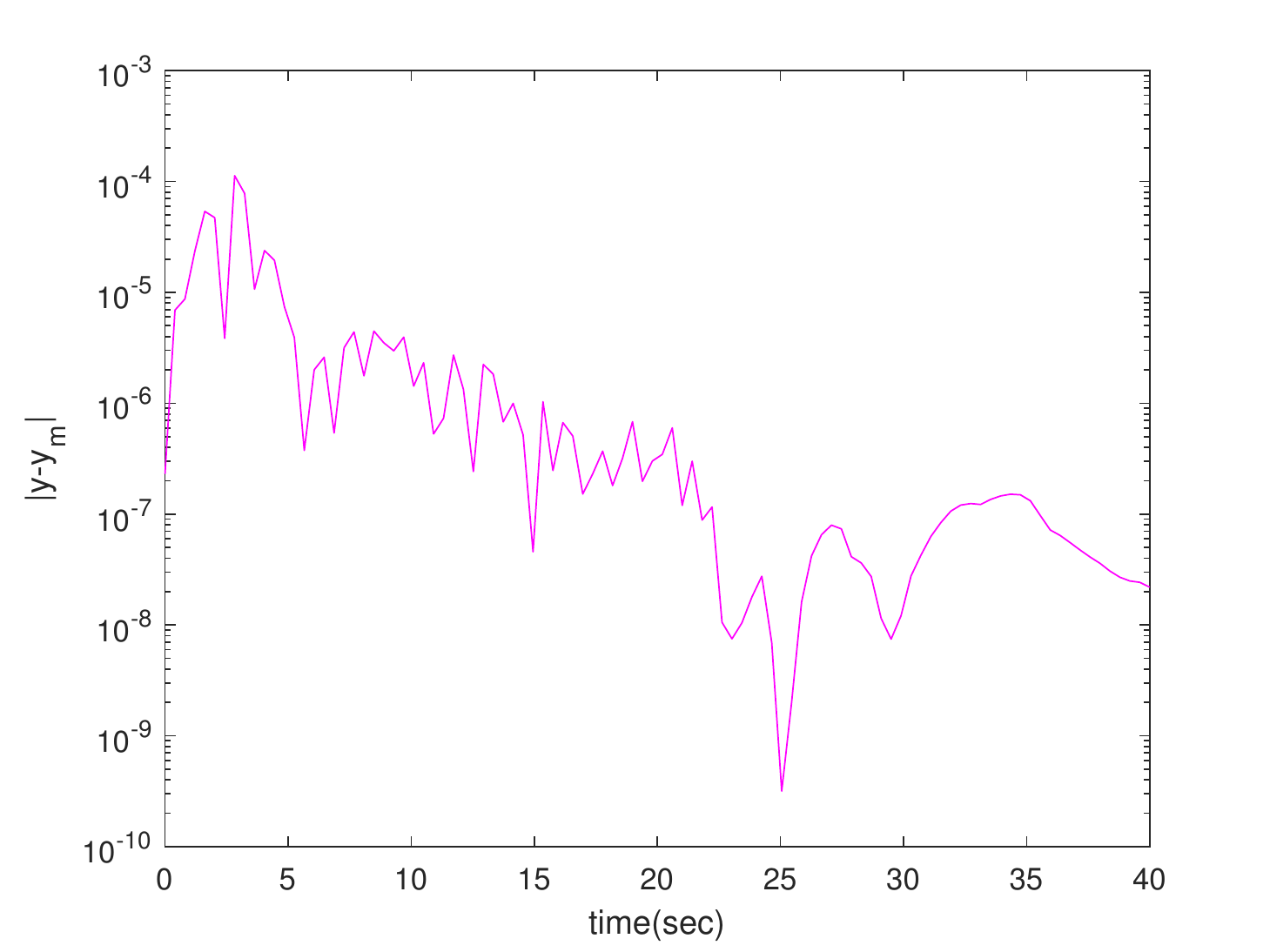}
	\caption{Left: time domain response for the stabilized system (input 1 to output 1) of original and reduced systems. Right: the error norm $\|y-y_m\|.$}
	\label{afstab500_1}
\end{figure}

\begin{figure}[H]
	\includegraphics[width=6.5cm,height=5.1cm]{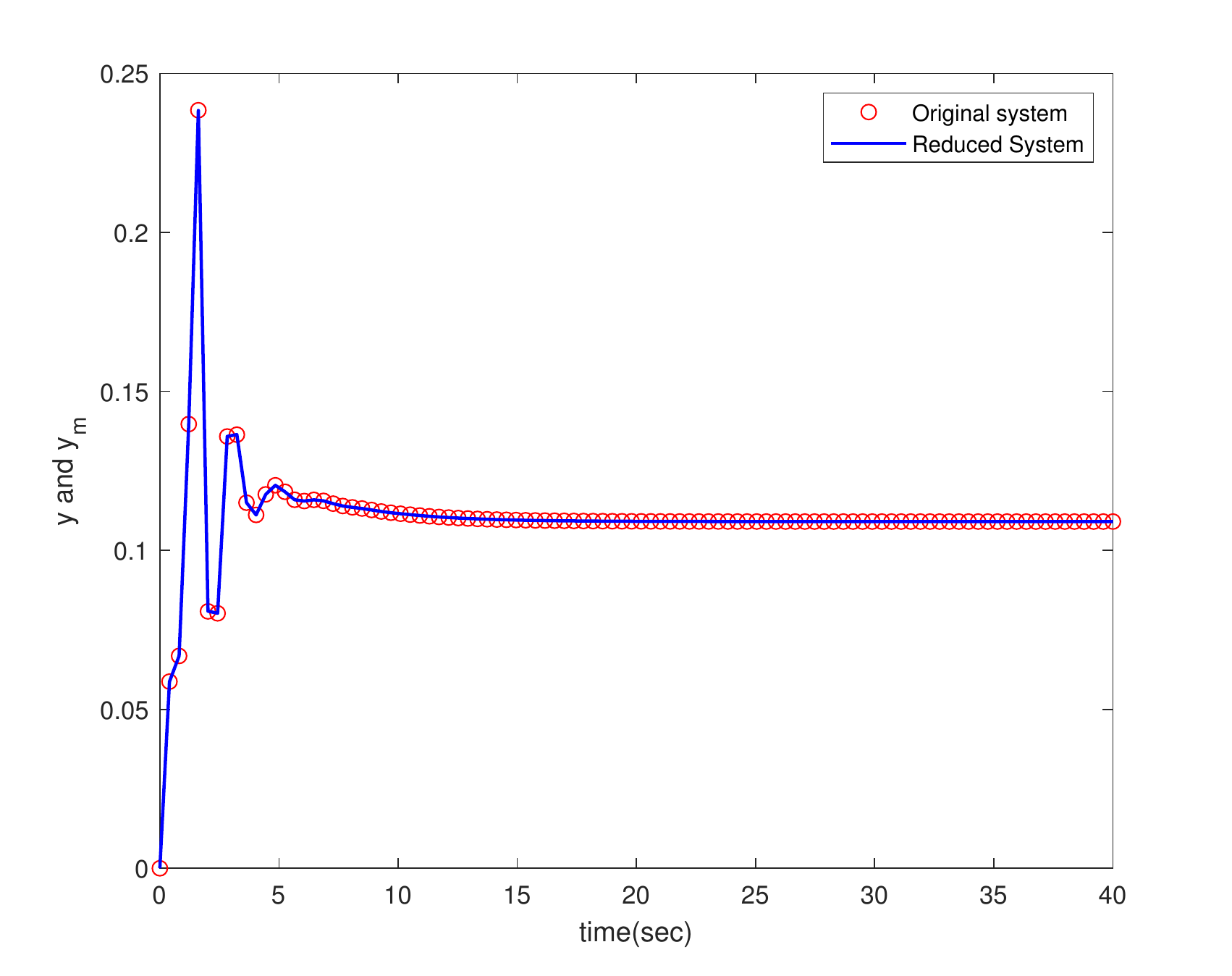}
	\hfill
	\includegraphics[width=6.5cm,height=5.1cm]{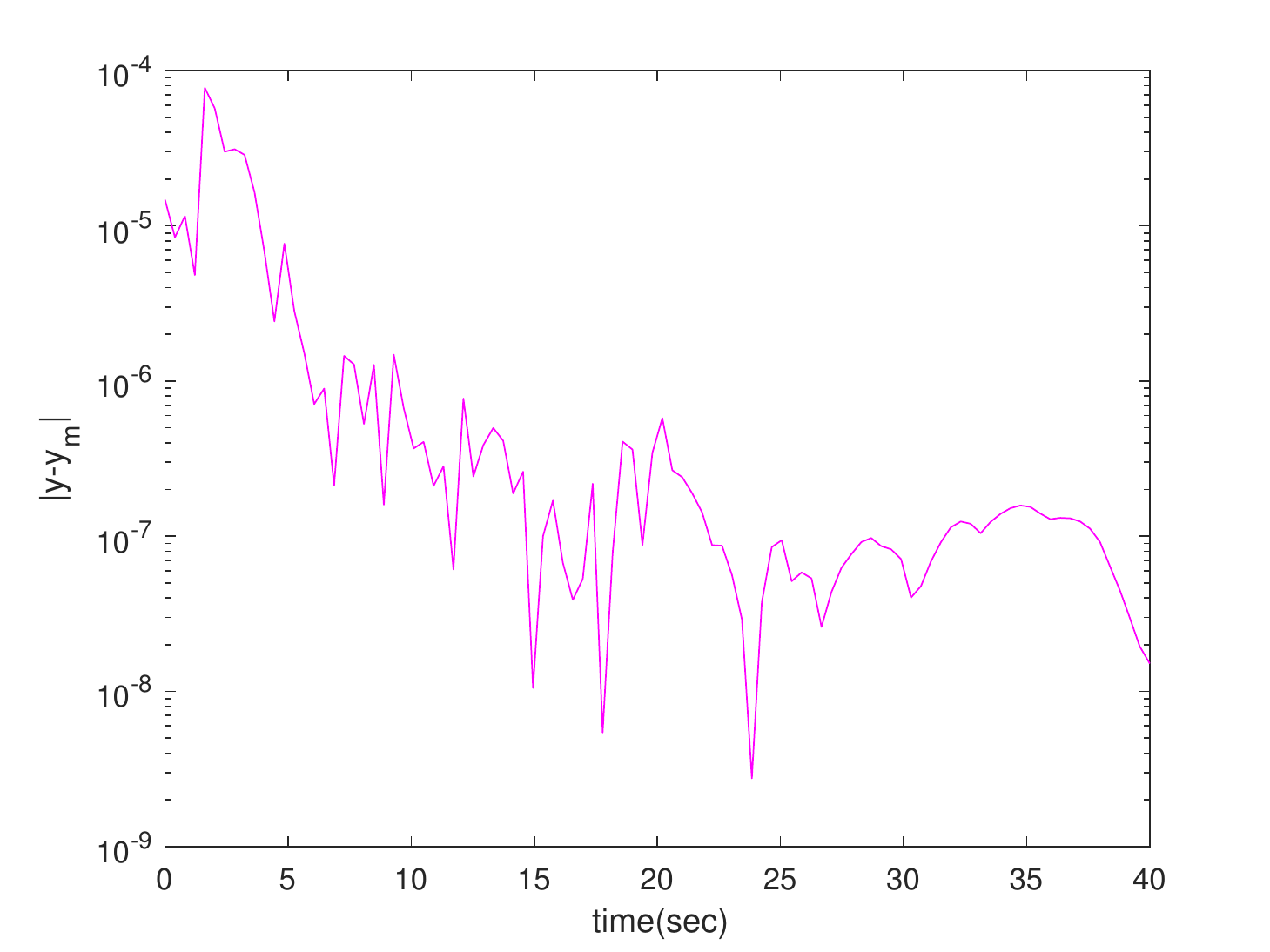}
	\caption{Left: time domain response for the stabilized system (input 2 to output 2) of original and reduced systems. Right: the error norm $\|y-y_m\|.$}
	\label{afstab500_2}
\end{figure}

We recall here the matrix feedback $K$ required to stabilize our original system (\ref{syscompact}). The control vector is  given by
$$\textbf{u}(t)=-K\bv(t) \quad \text{where} \quad K=B^T X_m M.$$
The matrix $X_m$ is the approximate solution to GARe (\ref{riccatinew}). We use the relation (\ref{ZZT}) that allows us to store $X_m$ in a efficient way, then the feedback matrix has the following form
$$K=B^T Z_m \, Z_m^T M.$$
The Reynolds number chosen here $\text{Re}=400$ and $500$ makes our original system (\ref{DAE1}) unstable as we mentioned earlier. We plug in the input ${\bf u}(t)$ in the unstable original system (\ref{DAE1}) to get the stabilized system described as follows 
\begin{subequations}
	\label{stabsys}
	\begin{align} 
		M \, \dfrac{d}{dt}\textbf{v}(t) &= (A-BK)\textbf{v}(t)+G\textbf{p}(t), \\ 
		0 &= G^T\textbf{v}(t), \\
		y&=C\, \textbf{v}(t).
	\end{align}
\end{subequations}
To show the effectiveness of the constructed feedback matrix $K$, we establish a time domain response simulation. In all examples (before and after stabilization) we use the same constant unit as input actuation. We use matrices corresponding to level 1 of discretization in Table \ref{tabnvnp} and we set $m=120$. For each $\text{Re}=400$ and $500$,  we first present the time domain response  of the original and reduced systems and  then we plot the time domain response associated with the stabilized system (\ref{stabsys}) and its reduced one. It is important to notice that while we perform the reduction process to the stabilized system (\ref{stabsys}),  using the extended block Krylov subspace method described in Section \ref{sec3}, we have to solve at each iteration the following saddle point problem
	\begin{align}
		\label{SDP}
		\underbrace{\begin{bmatrix}
				A-BK & G \\ 
				G^T & 0
		\end{bmatrix}}_{\bf \widehat{A}}  \begin{bmatrix}
			w \\ 
			\star
		\end{bmatrix} = \begin{bmatrix}
			z \\ 
			0
		\end{bmatrix},
	\end{align}
	Notice that the product $B K$ is dense and this in fact what makes the block  $(1,1)$ of $ \bf \widehat{A}$ dense too. To avoid this problem of density that can make our computation infeasible, we rewrite the saddle point problem (\ref{SDP}) in a low-rank form
	$$\left(\underbrace{\begin{bmatrix}
			A & G \\ 
			G^T & 0
	\end{bmatrix}}_{\bf A} -\underbrace{\begin{bmatrix}
			B \\ 
			0
	\end{bmatrix}}_{\bf B} \underbrace{\begin{bmatrix}
			K & 0
	\end{bmatrix}}_{\bf K}\right) \begin{bmatrix}
		w \\ 
		\star
	\end{bmatrix} = \begin{bmatrix}
		z \\ 
		0
	\end{bmatrix},$$
	and then we use the {\it Sherman-Morrison-Woodbury formula} \cite{golub} 
	$$({\bf A-BK})^{-1}= \bf A^{-1} + A^{-1}B(I-KA^{-1}B)^{-1} K\,A^{-1}.$$
	Besides solving the small dense matrix $\bf (I_{n_b}-KA^{-1}B)$ with right hand side $\bf K$ we need to solve $\bf A^{-1} B$ and $\bf A$ with the right hand side $[z,\,  0]^T$, and this can be done easily by adding the $n_b$ columns $\bf B$ to the matrix $[z,\,  0]^T$, and then instead of solving the problem (\ref{SDP}) with $\bf \widehat{A}$  one can solve the following saddle point problem
	$$\begin{bmatrix}
		A & G \\ 
		G^T & 0
	\end{bmatrix}\begin{bmatrix}
		w \\ 
		\star
	\end{bmatrix}=\begin{bmatrix}
		z & B\\ 
		0 & 0
	\end{bmatrix},$$
	using $"\backslash"$, a built-in MATLAB function.\\
In Figures \ref{bfstab400} and \ref{bfstab500}, we can see that for both cases $\text{Re}=400$ and $500$, the time domain simulation of the original and reduced systems show stability and a good accuracy of the reduced output compared to the original one. However, after  $t=30s$ some oscillations appear due to the instability of our original system. We can also see from the right parts of Figure \ref{bfstab400} and Figure  \ref{bfstab500}  that the error-norm   $\|y-y_m\|$ increases as the time  increases and this is due to the fact that our reduced system loses its accuracy caused by the instability that characterizes the original system. The performance of the matrix feedback  allows us to stabilize the unstable system. This is shown  in Figures \ref{afstab400_1}, \ref{afstab400_2}, \ref{afstab500_1} and \ref{afstab500_2} using two different Reynolds number $\text{Re}=400$ and $\text{Re}=500$. In the left side of these figures we display the time domain responses of the original and reduced stabilized systems of 1st input to 1st output in Figures \ref{afstab400_1} and \ref{afstab500_1} with $\text{Re}=400$ and $\text{Re}=500$ respectively, and also of 2nd input to 2nd output in Figures \ref{afstab400_2} and \ref{afstab500_2} with $\text{Re}=400$ and $\text{Re}=500$ respectively. One can notice that the figures illustrate a good accuracy of the reduced output compared to the original one. Moreover, it can be seen that after few oscillations that end in $t=5s$, the output of the stabilized system stabilize at constant values. On the right hand side of Figures \ref{afstab400_1}, \ref{afstab400_2}, \ref{afstab500_1} and \ref{afstab500_2}, we show  the error in the outputs for the same inputs and we notice that after the stabilization, the error $\|y-y_m\|$ does not increase as the time increases which was not the case before stabilization. This proves the accuracy of our method  of constructing a feedback matrix for stabilization using the EBARA Algorithm \ref{ERBALA}.

\section*{Conclusion}
Navier-stokes equations (NSEs) are considered as the pillars of fluid mechanics. A spatial discretization of the linearized NSEs around a steady state leads to a high dimension descriptor system of index-2 presented by a set of differential algebraic equations (DAEs). In this paper, we proposed a projection Krylov-based method to reduce this large dimension system. Our proposed method is based essentially on an extended block Arnoldi algorithm that allows us to build an efficient reduced system with a reasonable cost of computations. The system of NSEs lost its stability when  Reynolds numbers are large  and then  we need  stabilization techniques. One of the methods for stabilization that we used here is by solving an LQR problem based on a Riccati feedback approach. We suggested an extended Krylov-based method to solve the obtained  large-scale algebraic Riccati  equation and the obtained   numerical solution is the key to design a controller described by  a feedback matrix.

\end{document}